\documentclass[journal]{IEEEtran}%,doublespace
\usepackage{mathrsfs}
\input{epsf.sty}
\usepackage{epsfig}
\usepackage[dvipsnames,usenames]{color}
\usepackage{amsmath}
\usepackage{amssymb}
\usepackage{graphicx}
\usepackage{graphics,color}

\newtheorem{theorem}{Theorem}
\newtheorem{lemma}{Lemma}

\newtheorem{definition}{Definition}
\newtheorem{remark}{Remark}

\newcommand{\onetom}{1,\cdots,m}
\newcommand{\oneton}{1,\cdots,n}
\newcommand{\onetoK}{1,\cdots,K}

\begin{document}
\title {Cluster consensus in discrete-time networks of multi-agents with inter-cluster nonidentical inputs}
\author{Yujuan Han, Wenlian Lu,~\IEEEmembership{Member,~IEEE,}, Tianping Chen,~~\IEEEmembership{Senior~Member,~IEEE}\thanks{This work is jointly supported by the National Key Basic
Research and Development Program (No. 2010CB731403), the National Natural Sciences Foundation of China under Grant Nos. 61273211, 60974015, and 61273309, the Foundation for the Author of National Excellent Doctoral Dissertation of PR China No. 200921, the Marie Curie International Incoming Fellowship from the European Commission (FP7-PEOPLE-2011-IIF-302421), Shanghai Rising-Star Program (No. 11QA1400400), and also the Laboratory of Mathematics for Nonlinear Science, Fudan University. }\thanks{Yujuan Han is with School of Mathematical Sciences, Fudan University, No. 220 Handan Lu, Shanghai 200433, China.}
\thanks{Wenlian Lu is School of Mathematical Sciences, Fudan University, No. 220 Handan Lu, Shanghai 200433, China, and the Max Planck Institute in Mathematics in the Sciences, Inseltr. 22, Leipzig D04103, Germany.}
\thanks{Tianping Chen is with the School of Computer Sciences and the School of Mathematical Sciences, Fudan University, No. 220 Handan Lu, Shanghai 200433, China.}}
\date{}

\maketitle
\begin{abstract}
In this paper, cluster consensus of multi-agent systems is
studied via inter-cluster nonidentical inputs. Here, we consider general graph topologies, which might be time-varying.
The cluster consensus is defined by
two aspects: the intra-cluster synchronization, that the state
differences between each pair of agents in the same cluster converge to zero, and
inter-cluster separation, that the states of the agents in different
clusters are separated. For intra-cluster synchronization, the
concepts and theories of consensus including the spanning trees,
scramblingness, infinite stochastic matrix product and Hajnal
inequality, are extended. With them,  it is proved that if the graph
has cluster spanning trees and all vertices self-linked, then static
linear system can realize intra-cluster synchronization. For the
time-varying coupling cases, it is proved that if there exists $T>0$
such that the union graph across any $T$-length time interval has
cluster spanning trees and all graphs has all vertices self-linked,
then the time-varying linear system can also realize intra-cluster
synchronization. Under the assumption of common inter-cluster
influence, a sort of inter-cluster nonidentical inputs are utilized to realize
inter-cluster separation, that each agent in the same cluster
receives the same inputs and agents in different clusters have
different inputs. In addition, the boundedness of the infinite sum
of the inputs can guarantee the boundedness of the trajectory.
As an application, we employ a modified non-Bayesian social learning model
to illustrate the effectiveness of our results.
\end{abstract}
\begin{IEEEkeywords}
Cluster Consensus, Multi-agent System, Cooperative Control, Linear System, Non-Bayesian Social Learning
\end{IEEEkeywords}

\IEEEpeerreviewmaketitle
\section{Introduction}

In recent years, the multi-agent systems have broad
applications \cite{Vidal03,Cortes03,Fax04}. In particular, the consensus problems of multi-agent systems
have attracted increasing interests from many fields, such as
physics, control engineering, and biology \cite{Reynolds87,Vicsek95,Xiao04}. In network of
agents, \emph{consensus} means that all agents will converge to some
common state. A consensus algorithm is an interaction rule how
agents update their states. Recently, the consensus algorithm has also been used in social learning models.
Social learning focuses on the opinion dynamics in the society, which has attached a growing interests.
In social learning models, individuals engage in communication with their neighbors in order to learn from their experiences. For more details, we refer readers to see \cite{Acemoglu}-\cite{Liu qipeng}. A large amount of papers concerning consensus algorithms have been
published \cite{Ren05,Moreau05,Porfiri07,Hatano05,Liu08}, most of which focused on the average
principle,i.e., the current state of each agent is an average of the
previous states of its own and its neighbors, which is implemented by
communication between agents and can be described by the following difference equations for the discrete-time cases:
\begin{eqnarray}
x_{i}(t+1)=\sum_{j=1}^{n}A_{ij}x_{j}(t),~i=\oneton,\label{consensus}
\end{eqnarray}
where $x_{i}(t)$ denotes the state of agent $i$ and $A=[A_{ij}]_{i,j=1}^{n}$ is a {\em stochastic matrix}. For a survey, we refer readers to \cite{Olf07} and the references therein.

To realize consensus, the stability of the underlying dynamical system is curial. Since the network can be regarded as a graph, the issues can be depicted by the
graph theory. In the most existing literature, the concept of spanning tree is widely use to describe the communicability between agents in networks that can guarantee the consensus of (\ref{consensus}). See \cite{Ren04,Lu09,Liu09}.

It is widely known that the movement or/and defaults of the agents may lead the graph topology changing through time. So, it is inevitable to study the stability of the consensus algorithm in a time-varying environment, which can be described by the following time-varying linear system:
\begin{eqnarray}
x_{i}(t+1)=\sum_{j=1}^{n}A_{ij}(t)x_{j}(t),~i=\oneton,\label{TV-consensus}
\end{eqnarray}
where each $A(t)=[A_{ij}(t)]_{i,j=1}^{n}$ is a stochastic matrix. There were a lot of literature, in which the stability analysis of (\ref{TV-consensus}) are investigated. Most of their results can be derived from the theories of infinite nonnegative matrix product and ergodicity of inhomogeneous Markov chain. Among them, the followings should be highlighted. In \cite{Hajnal,Haj2}, the celebrated Hajnal's inequality was introduced and its generalized form was proposed in \cite{Shen}, to describe the compression of the differences among rows in a stochastic matrix when multiplied by another stochastic matrix that is scrambling. In \cite{Wolfowitz}, it was proved that a scrambling stochastic matrix could be obtained if a certain number of stochastic matrices that have spanning trees for their corresponding graphs were multiplied. So, in most of the papers involving stability analysis of (\ref{TV-consensus}), the sufficient conditions were expressed in terms of spanning trees in the union graph across time intervals of a given length. See \cite{Moreau05,Liu09} and the references therein. Besides, communication delays were also widely investigated \cite{Lu09,Olfati04,Liu10} and nonlinear consensus algorithms were proposed \cite{Liu108}.

All of the papers mentioned above concerns the complete consensus that the states of all agents converge to a common consistent state. However, this paper considered a more general
phenomenon, cluster consensus. This phenomenon
is observed when the agents in networks are
divided into several groups, called {\em clusters} in this paper by the way
that synchronization among the same cluster but the agents in different cluster have different
state trajectories. Cluster consensus (synchronization) is considered to be more
momentous in brain science \cite{schnitzler}, engineering control \cite{Passino}, ecological
science \cite{Montbrio}, communication, engineering \cite{Rulkov}, social
science \cite{Stone}, and distributed computation \cite{Hwang}.

In this paper, we define the cluster consensus as follows. Firstly, we divide the set of agents, denoted by $\mathcal V$, into disjoint
clusters, $\mathcal C=\{\mathcal C_{1},\cdots,\mathcal C_{K}\}$, with the properties:
\begin{enumerate}
\item $\mathcal C_{p}\bigcap\mathcal C_{q}=\emptyset$ for each $p\ne q$;
\item $\bigcup_{p=1}^{K}\mathcal C_{p}=\mathcal V$.
\end{enumerate}
Secondly, letting $x(t)=[x_{1}(t),\cdots,x_{n}(t)]^{\top}\in\mathbb R^{n}$ denote the state trajectory of all agents, of which $x_{i}(t)$ represents the state of $i\in\mathcal V$, we define {\em cluster consensus} via the following aspects:
\begin{enumerate}
\item $x(t)$ is bounded;

\item We say that $x(t)$ {\em intra-cluster synchronizes}
if $\lim_{t\rightarrow\infty}|x_{i}(t)-x_{i'}(t)|=0$ for all $i,i'\in \mathcal{C}_{p}$ and $p=1,\cdots,K$;

\item We say that $x(t)$ {\em inter-cluster separates} if $\limsup_{t\rightarrow\infty}|x_{i}(t)-x_{j}(t)|>0$ holds for each pair of $i\in\mathcal{C}_{k}$ and $j\in\mathcal{C}_{l}$ with $k\neq l$.
\end{enumerate}
We say that a system reaches {\em cluster consensus} if each solution $x(t)$ is bounded and satisfies the intra-cluster synchronization and inter-cluster separation, i.e., the items 1-3 are satisfied.

For this purpose, we introduce the following linear discrete-time
system with external inter-cluster nonidentical inputs:
\begin{eqnarray}\label{Mu-uniformCondition}
x_{i}(t+1)=\sum_{j=1}^{n}A_{ij}x_{j}(t)+I_{i}(t),~~
        i\in \mathcal{C}_{p} , p=1,\cdots,K,
\end{eqnarray}
where $A=[A_{ij}]_{i,j=1}^{n}$ is a $n\times n$ stochastic matrix,
$I_{i}(t)$ are external scalar inputs and they are different with respect to clusters,
which is used to realize inter-cluster separation. Also, we consider time-varying couplings that lead the following time-varying linear system with inputs:
\begin{eqnarray}
&&x_{i}(t+1)=\sum_{j=1}^{n}A_{ij}(t)x_{j}(t)+I_{i}(t),~~
        i\in \mathcal{C}_{p},\nonumber\\
        &&~~~~~~~~~~~~~~~~~~~~~~~~~~~~ p=1,\cdots,K.\label{Mu-uniformCondition1}
\end{eqnarray}

{\bf Related Works.}
Up till now, most papers in the literature mainly concern the global consensus. For instance, in \cite{Olf07,Liu09}, the (global) conseus was studied, especially for multi-agent system with time-varying topologies. There are essential differences between global consensus and the cluster consensus conidered the current paper, which means synchronization among the same cluster but the agents in different cluster have different state trajectories.
In some recent papers \cite{Chen11'}-\cite{LuChaos}, the authors addressed the cluster (group) consensus in networks with multi-agents and \cite{Chen11'} showed that (\ref{TV-consensus}) can reach cluster consensus if the graph topology is
fixed and strongly connected and the number of clusters is equal to the period of agents. For continuous-time network with fixed topology, \cite{Yu09} proved that under certain protocol, the multi-agent network can achieve group consensus by discussing the eigenvalues and eigenvectors of the Laplacian matrix. \cite{Yu09'} investigated group consensus in continuous-time network with switching topologies. However, all of these papers had a strong restriction in graph topologies and one important insight of cluster consensus: inter-cluster separation, has not been deeply investigated yet. Closely relating to this paper, the authors\rq{} previous work \cite{LuChaos} studied cluster synchronization of coupled nonlinear dynamical system and proposed several ideas, like intra-cluster synchronization  and configuration of graph topologies that cause cluster synchronization, which are shared in the current paper.

{\bf Our Contributions.} In this paper, we derive sufficient conditions for cluster consensus in the sense of both (\ref{Mu-uniformCondition}) and (\ref{Mu-uniformCondition1}). Different from the Lyapunov approach used in \cite{LuChaos}, in the current paper, we used the algebraic theory of product of infinite matrices and graph theory to derive the main results.  The enhancements in this paper, in comparison with the literature involved with (global) consensus like \cite{Olf07,Liu09} as well as the literature involved with cluster synchronization, like \cite{LuChaos}, are as follows. (1). We extended the concept of consensus to the cluster consensus as we mentioned above and the core concept of the algebraic graph theory, spanning tree, that means all nodes in the graph has a common root (a node can access all other nodes in the graph), to the cluster spanning tree, as defined in Definition \ref{cluster spanning tree}. (2). The main approach Hajnal inequality is extended to a cluster Hajnal inequality as Lemma 4. Accordingly, the concept of scramblingness is extended to cluster scramblingness as described in Definition \ref{cluster scrambling}. (3). We make efforts to prove inter-cluster separation, that the agents in different cluster do not converge to the same states, which is out of the scopes of the existing literature, like either \cite{Olf07,Liu09} or \cite{LuChaos}.

This paper is organized as follows. In section 2, we present
some graph definitions and give some notations required in this
paper. In section 3, we firstly investigate the cluster consensus
problem in discrete-time system with fixed topologies and present
the cluster consensus criterion. Then we promote the criterion to
the discrete-time system with switching topologies in section 4.
An application is given in section 5 to verify
the theoretical results. We conclude this paper in section 6.

\section{Preliminaries}

In this section, we firstly recall some necessary notations and
definitions that are related to graph and matrix theories and then
generalize them into the cluster sense. We also present several
lemmas which will be used later. For more details about the
definitions, notations and propositions about the graph and matrix,
we refer readers to textbooks \cite{God,Horn}.

For a matrix $A$, denote $A_{ij}$ the elements of $A$ on the $i$th row
and $j$th column. If the matrix $A$ is denoted as the result of an
expression, then we denote it by $[A]_{ij}$. $A^{\top}$ denotes the
transpose of $A$. For a
set $S$ with finite elements, $\#S$ denotes the number of elements
in $S$. $E$ denotes the identity matrix with a proper dimension.
${\bf 1}$ denotes the column vector with all components equal to $1$
with a proper dimension. $\rho(A)$ denotes the set of eigenvalues of a square matrix $A$. $\|z\|$ denotes a vector norm of a vector $z$
and $\|A\|$ denotes the matrix norm of $A$ induced by the vector norm $\|\cdot\|$.

An $n\times n$ matrix $A$ is called a {\em stochastic matrix} if
$A_{ij}\ge 0$ for all $i$, $j$, and $\sum_{j=1}^{n}A_{ij}=1$ for
$i=1,\cdots,n$. A stochastic matrix $A$ is called
{\em scrambling} if for any $i$ and $j$, there exists $k$ such that
both $A_{ik}$ and $A_{jk}$ are positive.

A directed graph $\mathcal{G}$ consists of a vertex set
$\mathcal{V}=\{\oneton\}$ and a directed
edge set $\mathcal{E}\subseteq
\mathcal{V}\times \mathcal{V}$, i.e., an
edge is an ordered pair of vertices in $\mathcal{V}$.
$\mathcal N_{i}$ denotes the neighborhood of the vertex $v_{i}$,
i.e. $\mathcal{N}_{i}=\{j\in\mathcal V:~(j,i)\in\mathcal
E\}$. A (directed) {\em path} of length $l$ from vertex
$v_{i}$ to $v_{j}$, denoted by
$(v_{r_{1}},v_{r_{2}},\cdots,v_{r_{l+1}})$, is a sequence of $l+1$
distinct vertices $v_{r_{1}},\cdots, v_{r_{l+1}}$ with
$v_{r_{1}}=v_{i}$ and $v_{r_{l+1}}=v_{j}$ such that
$(v_{r_{k}},v_{r_{k+1}})\in \mathcal{E}(\mathcal{G})$. The graph
$\mathcal{G}$ contains a {\em spanning (directed) tree} if there
exists a vertex $v_{i}$ such that for all the other vertices $v_{j}$
there's a directed path from $v_{i}$ to $v_{j}$, and $v_{i}$ is
called the {\em root} vertex. Corresponding to the matrix
scramblingness, we say that $\mathcal G$ is scrambling if for any
pair of vertices $v_{i}$ and $v_{j}$ there exists a common vertex $v_{k}$
such that $(v_{k},v_{i})\in\mathcal E$ and $(v_{k},v_{j})\in\mathcal
E$. We say that $\mathcal G$ has self-links if
$(v_{i},v_{i})\in\mathcal E$ for all $v_{i}\in\mathcal V$.

{\em Ergodicity coefficient}, $\mu(\cdot)$,
was proposed to measure the scramblingness of a stochastic matrix.
In \cite{Hajnal,Haj2}, the {\em Hajnal diameter} $\Delta(\cdot)$ was introduced to
measure the difference of the rows in a stochastic matrix, and
established his celebrated Hajnal's inequality
$\Delta(AB)\le(1-\mu(A))\Delta(B)$, which indicated that the Hajnal
diameter of stochastic matrix product $AB$ strictly decreases w.r.t.
$B$, if $A$ is scrambling, i.e., $\mu(A)<1$. The definitions of
$\mu(\cdot)$ and $\Delta(\cdot)$ can be found in \cite{Hajnal,Shen}.

An $n \times n$ nonnegative matrix $A$ can be
associated with a directed graph
$\mathcal{G}(A)=\{\mathcal{V},\mathcal{E}\}$ in such a way that
$(v_{j},v_{i})\in \mathcal{E}$ if and only if
$A_{ij}>0$. With this correspondence, we also say $A$ contains a
spanning tree if $\mathcal{G}(A)$ contains a spanning tree. On the
other hand, for a given graph $\mathcal G_{1}$, we denote
$\mathcal A(\mathcal G_{1})=\{A|\mathcal G(A)=\mathcal G_{1}\}$ the
subset of stochastic matrices $A$ such that $\mathcal G(A)=\mathcal
G_{1}$.

For an infinite stochastic matrix sequence $\{A(t)\}_{t=1}^{\infty}$
with the same dimension, we use the following simplified symbol for
a successive matrix product from $t$ to $s$ with $s>t$:
\begin{eqnarray*}
A_{t}^{s}\triangleq A(s)A(s-1)\cdots A(t).
\end{eqnarray*}
For a constant matrix $A$, we denote its $t$-th power by $A^{t}$. \cite{Wolfowitz} proved that if each stochastic matrix
$A(t)$ has spanning trees and self-links, then $A_{t}^{s}$ is
scrambling if $s-t>n-1$, where $n$ is the dimension of the matrix $A(t)$ \cite{Wu}.

In this paper, we consider cluster dynamics of networks. First of
all, for a graph $\mathcal G=(\mathcal V,\mathcal E)$, we define a
$\emph{clustering}$, $\mathcal{C}$, as a disjoint division of the
vertex set, namely, a sequence of subsets of $\mathcal{V}$,
$\mathcal{C}=\{\mathcal{C}_{1},\cdots,\mathcal{C}_{K}\}$, that
satisfies: (i) $\bigcup_{p=1}^{K}\mathcal{C}_{p}=\mathcal{V}$; (ii)
$\mathcal{C}_{k}\bigcap\mathcal{C}_{l}=\emptyset$, $k\neq l$. Thus,
we are able to extend the concepts of graph and matrix mentioned above
to those in the cluster case.

\begin{definition}\label{cluster spanning tree}
For a given clustering
$\mathcal{C}=\{\mathcal{C}_{1},\cdots,\mathcal{C}_{K}\}$, we say
that the graph $\mathcal G$ has {\em cluster-spanning-trees} with
respect to (w.r.t.) $\mathcal C$ if for each cluster
$\mathcal{C}_{p}$, $p=\onetoK$, there exists a vertex
$v_{p}\in\mathcal V$ such that there exist paths in $\mathcal G$
from $v_{p}$ to all vertices in $\mathcal{C}_{p}$. We denoted this
vertex $v_{p}$ as the root of the cluster $\mathcal C_{p}$.
\end{definition}

It should be emphasized that the root vertex of $\mathcal C_{p}$ and
the vertices of the paths from the root to the vertices in $\mathcal
C_{p}$ do not necessarily belong to $\mathcal C_{p}$. It can be seen
that the root vertex of a cluster does not necessarily same with the roots of
other clusters. Therefore, the definition of the
cluster-spanning-tree can be regarded as a generalization of that of
spanning tree we mentioned above.

\begin{definition}\label{cluster scrambling}
For a given clustering
$\mathcal{C}=\{\mathcal{C}_{1},\cdots,\mathcal{C}_{K}\}$, we say
that $G$ is {\em cluster-scrambling} (w.r.t. $\mathcal C$) if for
any pair of vertices $v_{p_{1}},v_{p_{2}}\in\mathcal C_{p}$,
there exists a vertex $v_{k}\in \mathcal V$, such that both
$(v_{k},v_{p_{1}})$ and $(v_{k},v_{p_{2}})$ belong to $\mathcal E$.
\end{definition}
Similarly, one can see that Definition \ref{cluster scrambling} is a
generalization of that of scramblingness we mentioned above. For a pair of vertices that are located in different clusters, they are not necessary to have a common linked vertex.

To measure the spanning-scramblingness, as a generalization from
those in Hanjnal \cite{Hajnal,Haj2}, we define the cluster ergodicity coefficient
(w.r.t the clustering $\mathcal C$) of a stochastic matrix $A$ as
\begin{eqnarray*}
\mu_{\mathcal C}(A)=\min_{p=1,\cdots,K}\min_{i,j\in\mathcal{C}_{p}}\sum_{k=1}^{n}{\min(A_{ik},A_{jk})}.
\end{eqnarray*}
It can be seen that $\mu_{\mathcal C}(A)\in [0,1]$ and $A$ is
cluster-scrambling (w.r.t. $\mathcal C$) if and only if
$\mu_{\mathcal C}(A)>0$.

According to the definition of cluster consensus, we extend the definition of Hajnal diameter \cite{Hajnal,Haj2,Shen} to the cluster case:
\begin{definition}
For a matrix $A$, which has row vectors $A_{1},A_{2},\cdots,A_{n}$ and a given clustering $\mathcal C$, we define the cluster Hajnal diameter as
\begin{eqnarray*}
\Delta_{\mathcal C}(A)=\max_{p=1,\cdots,K}\max_{i,j\in C_{p}}\|A_{i}-A_{j}\|
\end{eqnarray*}
for some norm $\|\cdot\|$.
\end{definition}
It can be seen that $\Delta_{\mathcal C}(x)\rightarrow 0$ is equivalent to the intra-cluster synchronization.

Similar to the results and the proof of Theorem 5.1 in \cite{Wu}, we can prove
that the product of $n-1$ $n-$dimensional stochastic matrices, all with cluster-spanning-trees, is cluster-scrambling.

\begin{lemma}\label{scrambling}
Suppose that each $A(t)$, $t=1,\cdots,n-1$ is an $n$-dimensional
stochastic matrix and has cluster-spanning-trees (w.r.t. $\mathcal
C$) and self-links. Then the product $A_{1}^{n-1}$ is cluster-scrambling (w.r.t. $\mathcal C$), i.e., $\mu_{\mathcal C}(A)>0$.
\end{lemma}

See the proof in Appendices.

In \cite{Olf07}, it has been proved that if a stochastic matrix $A$ has
spanning trees and all nodes self-linked, then the power matrix $A^{n}$
converges to $\mathbf{1}\alpha$ for some row vector
$\alpha\in\mathbb R^{n}$. Here, we conclude that the convergence can
hold even without the spanning tree condition as a direct consequence from \cite{Horn}.
\begin{lemma} \label{convergence} If a stochastic matrix $A$ has positive diagonal elements, then
$A^{n}$ is convergent exponentially.
\end{lemma}

\section{ Cluster consensus analysis of discrete-time network with static coupling matrix}

\subsection{Invariance of the cluster consensus subspace}
To seek sufficient conditions for cluster consensus, we firstly consider the situation that if the initial data $x(0)=[x_{1}(0),\cdots,x_{n}(0)]^{\top}$ has already had the cluster synchronizing structure, namely, $x_{i}(0)=x_{j}(0)$ for all $i,j\in \mathcal C_{p}$ with $p=\onetoK$, then the cluster synchronization should be kept ,i.e., $x_{i}(t)=x_{j}(t)$ for all $i,j\in \mathcal C_{p}$ with $p=\onetoK$ and $t\ge 0$. In other words, the following subspace in $\mathbb R^{n}$ w.r.t. the clustering $\mathcal C$:
\begin{eqnarray*}
\mathbb S_{\mathcal C}=\bigg\{&&x=[x_{1},\cdots,x_{n}]^{\top}\in\mathbb R^{n}:~x_{i}=x_{j}\\
&& {\rm~for~all~} i,j\in \mathcal C_{p}{\rm~with~} p=\onetoK \bigg\},
\end{eqnarray*}
named {\em cluster-consensus subspace}, is invariant through (\ref{Mu-uniformCondition}).

It should be emphasized that $I_{i}(t)$ are different with respect to clusters,
which is used to realize inter-cluster separation.
\begin{definition}
We say that the input $I(t)$ is {\em intra-cluster identical} if $I_{i}(t)=I_{j}(t)$ for all $i,j\in\mathcal C_{p}$ and all $p=\onetoK$ and the stochastic matrix $A$ has
{\em inter-cluster common influence} if for each pair of $p$ and $p'$, $\sum_{j\in C_{p'}}a_{ij}$ is identical w.r.t. all $i\in \mathcal C_{p}$, in other words, $\sum_{j\in C_{p'}}a_{ij}$ only depends on the cluster indices $p$ and $p'$ but is independent of the vertex $i\in C_{p}$.
\end{definition}

One can see that if two stochastic matrices $A$ and $B$ which have inter-cluster common influence w.r.t. the same clustering $\mathcal C$, so does the product $AB$. In the following, similar to what we did in \cite{LuChaos}, we have
\begin{lemma}\label{invariant}
If the input is intra-cluster identical and the matrix $A$ has inter-cluster common influence,
then the cluster-consensus subspace is invariant through (\ref{Mu-uniformCondition}).
\end{lemma}
{\em Proof.} From the condition, we define
\begin{eqnarray*}
\beta_{p,p'}\triangleq\sum_{j\in\mathcal C_{p'}}a_{ij}
\end{eqnarray*}
for any $i\in\mathcal C_{p}$ and $
I_{p}(t)\triangleq I_{i}(t)$
for any $i\in\mathcal C_{p}$.

Assuming that $x(t)\in \mathbb S_{\mathcal C}$, we are to prove $x(t+1)\in\mathbb S_{\mathcal C}$, too. For this purpose, let $
x_{p}(t)$ be the identical state of the cluster $p$ at time $t$. Thus,
for each $\mathcal C_{p}$ and an arbitrary vertex $i\in\mathcal C_{p}$,
\begin{eqnarray*}
x_{i}(t+1)&=&\sum_{p'=1}^{K}\sum_{j\in\mathcal C_{p'}}a_{ij}x_{j}(t)+I_{i}(t)\\
&=&\sum_{p=1}^{K}\beta_{p,p'}x_{p'}(t)+I_{p}(t),
\end{eqnarray*}
which is identical w.r.t. all $i\in\mathcal C_{p}$. By induction, this completes the proof.

\subsection{Intra-cluster synchronization}

We assume a special sort of intra-cluster identical input as follows:
\begin{eqnarray}
I_{i}(t)=\alpha_{p}u(t)\label{input}
\end{eqnarray}
where $u(t)$ is a scalar function and $\alpha_{1},\cdots,\alpha_{p}$ are different constants.

Similar to the Hanjnal inequality given in \cite{Hajnal,Haj2,Shen}, we can prove

\begin{lemma}\label{Hajnal} Suppose that stochastic matrices $A$ and $B$ having the same dimension and inter-cluster common influence,
then
\begin{eqnarray*}
\Delta_{\mathcal C}(AB)\le (1-\mu_{\mathcal C}(A))\Delta_{\mathcal C}(B).
\end{eqnarray*}
 \end{lemma}

{\em Proof.~} The idea of the proof is similar to that of the main result in \cite{Shen}. Let
\begin{eqnarray*}
B=\left[\begin{array}{c}B_{1}\\\vdots\\B_{n}\end{array}\right],~H=AB=\left[\begin{array}{c}H_{1}\\\vdots\\H_{n}\end{array}\right]
\end{eqnarray*}
with $B_{i}=[B_{i1},\cdots,B_{in}]$ and $H_{i}=\sum_{k}a_{ik}B_{k}$, denoted by $[H_{i1},\cdots,H_{in}]$, for all $i=\oneton$.

For any pair of indices $i$ and $j$ belonging to the same cluster $\mathcal C_{p_{0}}$, we have
\begin{eqnarray*}
H_{i}=\sum_{p=1}^{K}\sum_{k\in C_{p}}a_{ik}B_{k},~H_{j}=\sum_{p=1}^{K}\sum_{k\in C_{p}}a_{jk}B_{k}.
\end{eqnarray*}
Let $d_{k}=\min\{a_{ik},a_{jk}\}$. Define a set of index vector:
\begin{eqnarray*}
W=\{w=[w_{1},\cdots,w_{K}]:~w_{p}\in\mathcal C_{p},~p=\onetoK\}.
\end{eqnarray*}
For each $w\in W$, we define following convex combinations of $B_{1},\cdots,B_{n}$:
\begin{eqnarray*}
G_{w}=\sum_{p=1}^{K}\bigg[\sum_{k\in\mathcal C_{p},~k\ne w_{p}}d_{k}B_{k}+(\beta_{p,p_{0}}-\sum_{k\in\mathcal C_{p},~k\ne w_{p}}d_{k})B_{w_{p}}\bigg].
\end{eqnarray*}
It can be seen that both $H_{i}$ and $H_{j}$ are in  the convex hull of $G_{w}$ for all $w\in W$. Therefore,
\begin{eqnarray*}
\|H_{i}-H_{j}\|\le\max_{w,w'\in W}\|G_{w}-G_{w'}\|.
\end{eqnarray*}
Combining with
\begin{eqnarray*}
\|G_{w}-G_{w^{'}}\|&\le&\sum_{p=1}^{K}(\beta_{p,p_{0}}-\sum_{k\in\mathcal C_{p}}d_{k})\|B_{w_{p}}-B_{w^{'}_{p}}\|\\
&\le&(1-\mu_{\mathcal C}(A))\Delta(B).
\end{eqnarray*}
we have
\begin{eqnarray*}
\|H_{i}-H_{j}\|\le(1-\mu_{\mathcal C}(A))\Delta_{\mathcal C}(B).
\end{eqnarray*}
Therefore, $\Delta_{\mathcal C}(H)\le(1-\mu_{\mathcal C}(A))\Delta_{\mathcal C}(B)$,
which completes the proof due to the arbitrariness of $(i,j)\in\mathcal C_{p}$ and $p=\onetoK$.
\begin{remark}
Lemma \ref{Hajnal} indicates that if $A$ has inter-cluster common
influence, then the cluster-Hajnal diameter of $Ax$ decreases. In
addition, if $A$ is cluster-scrambling, $\Delta_{\mathcal C}(Ax)$ is
strictly less than $\Delta_{\mathcal C}(x)$.
\end{remark}

Based on the previous lemmas, we give the following result concerning
intra-cluster synchronization of (\ref{Mu-uniformCondition}).

\begin{theorem}\label{intra-cluster}
Suppose that both $u(t)$ and $\sum_{k=1}^{t}u(k)$
are bounded, $I(t)$ is defined by (\ref{input}), and $A$ is a stochastic
matrix with inter-cluster common influence, $A$ has cluster-spanning
trees and all positive diagonal elements. Then for any initial condition
$x(0)$, (\ref{Mu-uniformCondition}) is bounded and can
intra-cluster synchronize.
\end{theorem}

{\em Proof.} Let
$x(k)=[x_{1}(k),\cdots,x_{n}(k)]^{\top}$ be the solution of (\ref{Mu-uniformCondition}), then
\begin{eqnarray}
x(t+1)= A^{t+1}x(0)+\sum_{k=0}^{t}A^{t-k}I(k)
\end{eqnarray}
where $I(t)=\varsigma u(t)$ with
$\varsigma=[\varsigma_{1},\cdots,\varsigma_{n}]^{\top}$ and
\begin{eqnarray}
\varsigma_{i}=\alpha_{p},~i\in\mathcal C_{p}.\label{varsigma}
\end{eqnarray}
There is some $Y>0$ such that $|u(t)|\le Y$, $|\sum_{k=0}^{t}u(k)|\le Y$ hold for all $t\ge 0$.

By Lemma \ref{convergence}, we have $A^{t}=
A^{\infty}+\epsilon(t)$, where $\|\epsilon(t)\|_{\infty}\le
M\lambda^{t}$ for some $M>0$ and $\lambda\in (0,1)$. Therefore,
\begin{eqnarray*}
&&\|x(t+1)\|\leq\|A^{t+1}x(0)\|_{\infty}+\|\sum_{k=0}^{t}A^{t-k}\varsigma u(k)\|\\
&&\le\|x(0)\|+\|A^{\infty}\varsigma\| |\sum_{k=0}^{t}u(k)|+M\sum_{k=0}^{t}\lambda^{t-k}|u(k)|\\
&&\le\|x(0)\|+\|A^{\infty}\varsigma\| Y+MY\frac{1}{1-\lambda}.
\end{eqnarray*}
which implies the solution of system (\ref{Mu-uniformCondition}) is bounded.

By Lemma \ref{scrambling}, we can find an integer $N_{1}$ such that
for all $ m\geq N_{1}$, $A^{m}$ are cluster-scrambling. Denote
$\eta=1-\mu(A^{N_{1}})$. For any $t$, let $t=pN_{1}+l$ with some
$0\le l<p$. We have
\begin{eqnarray*}
\Delta_{\mathcal C}(A^{t+1})\le \eta^{p}\Delta_{\mathcal C}(E_{n})
\end{eqnarray*}
which converges to zero as $t\to\infty$. In addition, since $A^{l}$
has inter-cluster common influence and $\Delta_{\mathcal
C}(\varsigma)=0$, then $\Delta_{\mathcal C}(A^{l}\varsigma)=0$ for all
$l\ge 0$ can be concluded. Therefore, we have $\Delta_{\mathcal
C}(x(t+1))\le \Delta_{\mathcal C}(A^{t+1}x(0))$ converges to zero as
$t\to\infty$. This completes the proof.

\subsection{Inter-cluster separation}

Under the conditions of Theorem \ref{intra-cluster}, the system can
intra-cluster synchronize, namely, the states within the same
cluster approach together. However, it is not known if the states in
different clusters will approach together, too. A simple
counter-example is that the matrix $A$ with the inter-cluster common
influence has (global) spanning trees with all diagonal elements positive and the
inputs  $\varsigma u(t)$ satisfies $\sum_{k=1}^{\infty}|u(k)|$
converges. In this case, $u(t)$ converges to zero and the
influence of the input to the system disappears. One
can see that $x(t)$ reaches a global consensus, i.e.,
$\lim_{t\to\infty}x(t)={\bf 1}\alpha$ for some scalar $\alpha$.

In this section, we investigate this problem by assuming that $u(t)$ is periodic with a
period $T$ and $\sum_{k=1}^{T}I(k)=0$, which guarantees that the sum
of $u(t)$ is bounded. Construct a new matrix: $B=[\beta_{p,q}]_{p,q=1}^{K}$, where
\begin{eqnarray}
\beta_{p,q}=\sum_{j\in\mathcal C_{q}}a_{ij},~~ i\in \mathcal
C_{p}\label{B}
\end{eqnarray}
It can be seen that $\beta_{p,q}$ is independent of the selection of
$i\in\mathcal C_{p}$.

Furthermore, we use the concept of ``genericality'' from the structural control theory \cite{Lin,Rein,Dion} to investigate
the inter-cluster separation. We define a set $\mathcal T(\mathcal
C,\mathcal G)$ w.r.t. a clustering $\mathcal C$ and a graph
$\mathcal G$, of which each element has form:
$\{B,\tilde{\varsigma},[u_{1},\cdots,u_{T-1}]\}$, where $B$ is defined in (\ref{B}) corresponding to the graph topology $\mathcal G$,
$\tilde{\varsigma}\in\mathbb R^{K}$ is the vector to identify each cluster and defined as:
\begin{eqnarray}
\tilde{\varsigma}_{p}=\alpha_{p},~p=\onetoK,\label{tvar}
\end{eqnarray}
and
$[u_{1},u_{2},\cdots,u_{T-1}]\in\mathbb R^{T-1}$ such that
\begin{eqnarray}
&&u(\theta+kT)=u_{\theta},~\theta=1,\cdots,T-1,\nonumber\\
&&{\rm and}~u(kT)=-\sum_{j=1}^{T-1}u_{j},~\forall~k\ge 0.\label{u}
\end{eqnarray}
We can
rewrite the system (\ref{Mu-uniformCondition}) as the following
compact form:
\begin{eqnarray}
x(t+1)=Ax(t)+\varsigma u(t),\label{compact}
\end{eqnarray}

\begin{definition}
We say that for a given set $\mathcal T(\mathcal C,\mathcal G)$ as
defined above, (\ref{compact}) is {\em generically} inter-cluster
separative (or cluster consensus) if for almost every triple
$\{B,\tilde{\varsigma},[u_{1},\cdots,u_{T-1}]\}\in\mathcal T(\mathcal C,\mathcal G)$ and almost
all initial $x(0)\in\mathbb R^{n}$, (\ref{compact}) can inter-cluster
separate (or cluster consensus).
\end{definition}

Before presenting a sufficient condition for generical inter-cluster
separation, we give the following simple lemma.

\begin{lemma}\label{common link}
Suppose that the stochastic matrix $A$ has the
inter-cluster common influence. Then, for any pair of cluster
$\mathcal C_{1}$ and $\mathcal C_{2}$, either there are no links
from $\mathcal C_{2}$ to $\mathcal C_{1}$; or for each vertex
$v\in\mathcal C_{1}$, there is at least one link from $\mathcal
C_{2}$ to $v$.
\end{lemma}

\begin{theorem}\label{inter-cluster}
Suppose that
\begin{enumerate}
\item every vertex in $\mathcal G$ has a self-link;

\item $\mathcal G$ satisfies the condition in Lemma \ref{common link}
w.r.t $\mathcal C$;

\item $\mathcal G$ has cluster-spanning-trees.
\end{enumerate}
Then (\ref{compact}) reaches
cluster consensus generically with respect to the set $\mathcal
T(\mathcal C,\mathcal G)$. In addition, the limiting consensus
trajectories are periodic, that is, there exist some scalar
periodic trajectories $v_{p}(t)$ with the period $T$ for each cluster $\mathcal C_{p}$, $p=\onetoK$,
such that $\lim_{t\to\infty}|x_{j}(t)-v_{p}(t)|=0$ if $j\in\mathcal
C_{p}$.
\end{theorem}

{\em Proof.} We firstly prove the asymptotic periodicity.
Recall
\begin{eqnarray}
x(t+1)= A^{t+1}x(0)+\sum_{k=0}^{t}A^{k}\varsigma u(t-k).
\end{eqnarray}

By Lemma \ref{convergence}, one can see that $A^{n}$ exponentially converges to $A^{\infty}$. Thus, we can
find $M>0$ and $\lambda\in (0,1)$ such that
$\|A^{t}-A^{\infty}\|\le M\lambda^{t}$. Let
$Y=\max_{k=1,\cdots,T}|u(k)|$. Thus, we have
\begin{eqnarray*}
&&\|x(t+lT+1)-x(t+1)\|\leq
\|(A^{t+lT+1}-A^{t+1})x(0)\|\\
&&+\|\sum_{k=0}^{t}A^{k}\varsigma[u(t+lT-k)-u(t-k)]\|\\
&&+\|\sum_{k=t+1}^{t+lT}A^{k}\varsigma u(t+lT-k)\|\\
&&=\|(A^{t+lT+1}-A^{t+1})x(0)\|
+\|\sum_{k=t+1}^{t+lT}A^{\infty}\varsigma u(t-k)\|\\
&&+\|\sum_{k=t+1}^{t+lT}[A^{k}-A^{\infty}]\varsigma u(t-k)\|\\
&&\leq
2M\lambda^{t}\|x(0)\|+MY\|\varsigma\|\sum_{k=t}^{t+lT}\lambda^{k}\\
&&=\bigg[2M\|x(0)\|
+MY\|\varsigma\|\frac{1}{1-\lambda}\bigg]\lambda^{t}
\end{eqnarray*}
for all $l$. Letting $t=mT+\theta-1$ for any $m\in\mathbb N$ and $\theta=1,\cdots,T$, we have
\begin{eqnarray*}
\|x((m+l)T+\theta)-x(mT+\theta)\|\le M_{1}\lambda^{mT}
\end{eqnarray*}
for some $M_{1}>0$. According to the Cauchy convergence principle, each $x(\theta+kT)$, $\theta=1,\cdots,T$, converges to some value as $k\to\infty$ exponentially, which implies that there exist T-periodic functions $v_{p}(t)$,
$p=1,\cdots,K$, such that
$|x_{j}(t)-v_{p}(t)|\rightarrow 0$ exponentially, if $j\in\mathcal
C_{p}$.

Now, we will prove  the consensus states in different clusters are different generically.

Since each cluster synchronizes, we can pick a single vertex state
from each cluster to represent the whole state of this cluster. We can divide the space $\mathbb R^{n}$ into the
direct sum of two subspaces: $\mathbb R^{n}=V_{1}\oplus V_{2}$, where $V_{1}$ denotes the right eigenspace of $A$ corresponding to the eigenvalue $1$ and $V_{2}$ denotes the right eigenspace of $A$ corresponding to all other eigenvalues. Since all diagonal elements of $A$ are positive, then the direct sum works and $AV_{i}\subset V_{i}$ holds for $i=1,2$. In addition, since the column vectors in $A^{\infty}$ belong to $\mathcal S_{\mathcal C}$, $V_{1}\subset\mathcal S_{\mathcal C}$. So, $\mathbb R^{n}=\mathcal S_{\mathcal C}+V_{2}$.

For any initial data $x(0)\in\mathbb R^{n}$, we can find $y^{0}\in\mathcal S_{\mathcal C}$ with the decomposition $x(0)=y^{0}+x(0)-y^{0}$ such that $x(0)-y^{0}\in V$. Consider the following system restricted to $\mathcal S_{\mathcal
C}$:
\begin{eqnarray*}
y(t+1)=Ay(t)+I(t),~y(0)=y_{0}.
\end{eqnarray*}
where $y(t)\in\mathcal S_{\mathcal C}$ for all $t$.

Let $\delta x(0)=x(0)-y^{0}\in V_{2}$. We have
\begin{eqnarray*}
\delta x(t+1)=A\delta x(t), \delta{x(0)}=x(0)-y^{0},
\end{eqnarray*}
which implies that $\lim_{t\rightarrow \infty}\delta x(t)=0$, that is, $\lim_{t\rightarrow \infty}x(t)=\lim_{t\rightarrow \infty}y(t)$. Therefore, we only need to discuss $y(t)\in \mathcal{S}_{\mathcal{C}}$.

Since each component of $y(t)$ in the same cluster is identical, we can
pick a single component from each cluster to lower-dimensional
column vector $\tilde{y}\in\mathbb R^{K}$ with $\tilde{y}_{p}=y_{i}$
for some $i\in\mathcal C_{p}$. Because of the inter-cluster common influence condition, we have
\begin{eqnarray}
\tilde{y}(t+1)=B\tilde{y}(t)+\tilde{\varsigma} u(t)\label{sum}
\end{eqnarray}
where $B$ is defined in (\ref{B}) and
$\tilde{\varsigma}=[\alpha_{1},\cdots,\alpha_{K}]^{\top}$. The
inter-cluster separation is equivalent to investigate the
separation among components of $\tilde{y}(t)$. One can see
that for almost every $B$, $B$ has $K$ distinguishing left
eigenvectors, denoted by $\phi_{1},\cdots,\phi_{K}$, corresponding
to eigenvalues $\nu_{1},\cdots,\nu_{K}$ (possibly overlapping). So,
for almost every $B$ with $K$
left eigenvectors, let us write down the solution
(\ref{sum}) at time $nT$ as follows:
\begin{eqnarray*}
\tilde{y}(nT+1)&=&B^{nT+1}\tilde{y}(0)+\sum_{k=0}^{nT}B^{nT-k}\tilde{\varsigma}u(k)\\
&\rightarrow&Z_{1}\tilde{y}(0)+Z_{2}\tilde{\varsigma},~{as~}n\to\infty,
\end{eqnarray*}
where
\begin{eqnarray*}
Z_{1}=\lim_{n\to\infty}B^{nT+1},~Z_{2}=\lim_{n\to\infty}\bigg[\sum_{k=0}^{nT}B^{nT-k}u(k)\bigg].
\end{eqnarray*}
From Lemma \ref{convergence},
$Z_{1}$ does exist. Combined with $\sum_{k=0}^{nT-1}u(k)=0$, we can conclude that $Z_{2}$ exists, too.

For an arbitrary fixed pair of $(p,q)$, with $p, q=\onetoK$ and $p\ne q$, we are to show $Z_{2}$ can generically
have different $p$-th and $q$-th components. In fact, for each $k_{1}$ with
$|\nu_{k_{1}}|<1$, noting that its associated left eigenvector is $\phi_{k_{1}}$, we have
\begin{eqnarray*}
&&\phi_{k_{1}}\sum_{k=0}^{nT}B^{nT-k}u(k)\\
&&=\phi_{k_{1}}\sum_{k=0}^{T-1}u(k)\sum_{z=0}^{n-1}
\nu_{k_{1}}^{zT+k}+\phi_{k_{1}}\nu_{k_{1}}^{nT}u(0)\\
&&\rightarrow\phi_{k_{1}}\frac{\sum_{k=0}^{T-1}u(k)\nu_{k_{1}}^{k}}{1-\nu_{k_{1}}^{T}},~{\rm
as}~n\to\infty.
\end{eqnarray*}
For each $k_{2}$ with $|\nu_{k_{2}}|=1$, noting its associated left-eigenvector is $\phi_{k_{2}}$, according to the fact that all
diagonal elements in $B$ are positive, from \cite{Horn}, we have $\nu_{k_{2}}=1$ indeed.
Then, we have
\begin{eqnarray*}
\phi_{k_{2}}\sum_{k=0}^{nT}B^{nT-k}u(k)=\phi_{k_{2}}\sum_{k=0}^{nT}u(k)=\phi_{k_{2}}u(0)=\phi_{k_{2}}u_{T}.
\end{eqnarray*}
So, for almost $[u_{1},\cdots,u_{T-1}]\in\mathbb R^{T-1}$, the eigenvectors of $Z_{2}$ are the same with $B$ and the corresponding eigenvalues are $u_{T}$ and $\sum_{k=0}^{T-1}u(k)\nu_{p}^{k}/(1-\nu_{p}^{T})$. For almost every realization of $[u_{i}]_{i=1}^{T-1}$ and $B$, none of them is zero, which implies that $Z_{2}$ is nonsingular. That means it is impossible for each pair of its rows to be identical. So, for almost all $\tilde{\varsigma}$, the $p$-th and $q$-th component of $Z_{2}$ are not identical. Equivalently, for almost every $\tilde{\zeta}$, $Z_{2}\tilde{\zeta}$ has no pair of components identical. Therefore, we conclude that for almost every $x_{0}$, associated with almost every $\tilde{y}(0)$, each pair of components in $Z_{1}\tilde{y}(0)+Z_{2}\tilde{\zeta}$ are not identical.

We can arbitrarily select the cluster pair $(p,q)$ and the exception cases of the statements above are within a set of $\mathcal T(\mathcal G,\mathcal C)$ with Lebesgue measure zero. Since any finite union of sets with Lebesgue measure zeros still has Lebesgue measure zero, we conclude that $\lim_{n\to\infty}\tilde{y}(nT+1)$ has no identical components generically, which implies that the states of any two clusters in $\lim_{n\to\infty}y(nT+1)$ are not identical generically. This completes the proof.

\begin{remark}
In the current paper, we make efforts to prove the inter-cluster separation rigoroursly; however, in \cite{LuChaos}, the inter-cluster separation was not touched (but only assumed). We argue that for general nonlinear coupled system (models in \cite{LuChaos}), proving the inter-cluster separation is very difficult , if it was not impossible.
\end{remark}

\section{Cluster-consensus in discrete-time network with switching topologies}

In this section, we study the cluster-consensus in network with switching topologies described as the following time-varying linear system:
\begin{eqnarray}
&&x_{i}(t+1)=\sum_{j=1}^{N}A_{ij}(t)x_{j}(t)+I_{i}(t)~~
    \forall i\in \mathcal{C}_{p},\nonumber\\
&&~~~~~~~~~~~~~~~~~~~~~~~~~~~~~~p=\onetoK,\label{time varying}
\end{eqnarray}
where $A(t)$ is associated with a graph from the graph set $\Upsilon=\{\mathcal G_{1},\cdots,\mathcal G_{m}\}$ w.r.t. a given clustering $\mathcal C$, each of which satisfy the property ${\bf \mathcal A}$: for each pair $p$ and $q$ of cluster indices,
\begin{enumerate}
\item there are no links from $\mathcal C_{q}$ to $\mathcal C_{p}$ in each graph $\mathcal G_{l}$, $l=\onetom$,
\item or for each vertex $v\in\mathcal C_{p}$ and each graph $\mathcal G_{l}$, $l=\onetom$, there is at least one link from $\mathcal C_{q}$ to it.
\end{enumerate}
For the matrix sequence $A(t)$, we have the following assumptions:
\begin{itemize}
\item ${\bf \mathcal B}_{1}$: There is a positive constant $e>0$ such that for each pair $i,j$ and $t$, either $A_{ij}(t)=0$ or $A_{ij}\ge e$ holds;
\item ${\bf \mathcal B}_{2}$: $A_{ii}(t)\ge e$ holds for all $i=\oneton$ and $t\ge 0$;
\item ${\bf \mathcal B}_{3}$ ({\em inter-cluster common influence}): There exists a $\mathbb R^{n,n}$ stochastic matrix $B(t)=[b_{p,q}(t)]_{p,q=1}^{K}$, possibly depending on time, such that
\begin{eqnarray}
\sum_{j\in C_{q}}A_{ij}(t)=b_{p,q}(t)\label{B1}
\end{eqnarray}
holds for all $i\in\mathcal C_{p}$ and each $p,q=\onetoK$;
\item ${\bf \mathcal B}_{3}^{*}$ ({\em static inter-cluster common influence}): There exists a constant $\mathbb R^{n,n}$ stochastic matrix $B=[b_{p,q}]_{p,q=1}^{K}$, such that
\begin{eqnarray}
\sum_{j\in C_{q}}A_{ij}(t)=b_{p,q}\label{B1a}
\end{eqnarray}
holds for all $i\in\mathcal C_{p}$ and each $p,q=\onetoK$.
\end{itemize}
In other words, we define a graph set containing all possible graph induced by the matrix sequence $A(t)$. The graph set satisfies the property in Lemma \ref{common link} uniformly and each graph in the set either never occurs in the corresponding graph sequence induced by $A(t)$ or occurs frequently.

Then, we are in the position to give a sufficient condition for the cluster synchronization.

\begin{theorem}\label{tv intra-cluster}
Suppose that $\mathcal A$, $\mathcal B_{1}$, $\mathcal B_{2}$ and $\mathcal B_{3}$ hold. If
there exists an integer $L>0$ such that for any $L$-length time interval $[t,t+L)$, the union graph $\mathcal G[\sum_{i=t}^{t+L-1}A(i)]$ has cluster-spanning-trees, then the system (\ref{time varying}) cluster synchronizes.
\end{theorem}
{\em Proof.}
The solution of (\ref{time varying}) is
\begin{eqnarray*}
x(t+1)=A(t)x(t)+\varsigma u(t)
=A_{0}^{k}x(0)+\sum_{k=0}^{t} A_{k+1}^{t} u(t)\varsigma.
\end{eqnarray*}
Noting that the diagonal elements of each $A(t)$ are positive, we can see that the graph $\mathcal G(A_{t}^{t+L-1})$ contains all links in the union graph $\mathcal G(\sum_{k=t}^{t+L-1}A(k))$ and hence has cluster-spanning-trees and positive diagonal elements for all $t$. By Lemma \ref{scrambling}, we can conclude that there is an integer $N$ such that the graph $\mathcal G(A_{t}^{t+NL-1})$ is scrambling for all $t\ge 0$. Since the nonzero elements in each $A(t)$ is greater than some constant $e>0$, there is some $\delta>0$ such that
\begin{eqnarray*}
\inf_{t}\mu_{\mathcal C}(A_{t}^{t+NL-1})\ge\delta.
\end{eqnarray*}
Hence, for each $t$, we have
\begin{eqnarray*}
\Delta_{\mathcal C}(A_{0}^{t}x(0))\le(1-\delta)^{\lfloor\frac{t}{NL}\rfloor}\Delta_{\mathcal C}(x(0)),
\end{eqnarray*}
which converges to zero as $t\to\infty$. Here $\lfloor\cdot\rfloor$ denotes the floor function. Therefore, $\lim_{t\to\infty}\Delta_{\mathcal C}(A_{0}^{t}x(0))=0$.

Combining with the fact that $\Delta_{\mathcal C}(A^{s}_{t}\varsigma)=0$ holds for all $s\ge t$ and $\varsigma$,
we can conclude that the system (\ref{time varying}) intra-cluster synchronizes.

\begin{remark}
Due to the difference of the techniques used in \cite{LuChaos} and the current paper, the result of Theorem \ref{tv intra-cluster}  is impossible to extend to general coupled nonlinear system, as the models in \cite{LuChaos}, because a Lyapunov function for time-varying coupled systems is in general unable to be found.
\end{remark}

The inter-cluster separation can be derived by the same fashion of Theorem \ref{inter-cluster}.
\begin{theorem}
Suppose that ${\bf \mathcal A}$,  ${\mathcal B}_{1}$, $\mathcal B_{2}$ and ${\mathcal B}_{3}^{*}$ hold. If there exists an integer $L>0$ such that for any $L$-length time interval $[t,t+L)$, the union graph $\mathcal G[\sum_{i=t}^{t+L-1}A(t)]$ has cluster spanning trees. If the input $u(t)$ and $\sum_{k=0}^{t} u(k)$ are both bounded, then for any initial data $x(0)$, the solution of (\ref{time varying}) is bounded. In addition, if the input $u(t)$ is periodic with a period $T$ and satisfies $\sum_{k=1}^{T-1}u(k)=0$, (\ref{time varying}) reaches cluster consensus generically and each trajectory converges to a $T$-periodic one.
\end{theorem}
{\em Proof.} To prove the boundedness, we are to find a solution of (\ref{time varying}) that stays at $\mathcal S_{C}$ and is the limiting of $x(t)$. Similar to the proof of Theorem \ref{inter-cluster}, we can represent the limiting trajectory by a lower-dimensional linear system (\ref{sum}). The ${\bf\mathcal B}_{3}^{*}$ implies that this linear lower-dimensional system is static. So, we can prove its boundedness by the same way of the proof of Theorem \ref{intra-cluster}.

Define the Lyapunov exponent of the matrix sequence $A(t)$ as follows:
\begin{eqnarray*}
\lambda(v)=\overline{\lim}_{t\to\infty}\frac{1}{t}\log\bigg(\|A_{0}^{t}v\|\bigg).
\end{eqnarray*}
From the Pesin's theory \cite{Pesin}, the Lyapunov exponents can only pick finite values and provide a splitting of $\mathbb R^{n}$. Namely, there is a subspace direct-sum division:
\begin{eqnarray*}
\mathbb R^{n}=\oplus_{j=1}^{J}V_{j},
\end{eqnarray*}
and $\lambda_{1}>\cdots>\lambda_{J}$, possibly $J<n$, such that for each $v\in V_{j}$, $\lambda(v)=\lambda_{j}$. It can be seen that $\lambda_{1}=0$ since $A(t)$, $t\ge 0$, are all stochastic matrices.  Let $V=\oplus_{j>1}V_{j}$. We claim

{\em Claim 1}: $\mathbb R^{n}=\mathcal S_{\mathcal C}+V$.

We prove this claim in appendix. Therefore, for any $x(0)\in\mathbb R^{n}$, we can find a vector $y_{0}\in\mathcal S_{\mathcal C}$ such that $x(0)-y_{0}\in V$. Define a linear system:
\begin{eqnarray}
y(t+1)=A(t)y(t)+\varsigma u(t),~y(0)=y_{0}.\label{y}
\end{eqnarray}
Then, letting $\delta x(t)=x(t)-y(t)$, it should satisfy:
\begin{eqnarray*}
\delta x(t+1)=A(t)\delta x(t),~\delta x(0)=y(0)-x(0)\in V.
\end{eqnarray*}
Since $\delta x(0)\in V$, $\lambda(\delta x(0))<0$. This implies $\lim_{t\to\infty}
\delta x(t)=0$. So, $\lim_{t\to\infty}[x(t)-y(t)]=0$. We can rewrite the equation (\ref{y}) as a lower-dimensional linear system:
\begin{eqnarray}
\tilde{y}(t+1)=B\tilde{y}(t)+\tilde{\varsigma}u(t),\label{sum1}
\end{eqnarray}
which is same with (\ref{sum}). The  ${\bf\mathcal B}_{3}^{*}$ guarantees that the matrix $B$ is static.
So, the proof of boundedness of $\tilde{y}(t)$ is an overlap of that of Theorem \ref{intra-cluster}.

In addition, since $B$ is static, then the inter-cluster separation can be proved as an overlap of that of Theorem \ref{intra-cluster}.
Therefore, we can conclude that $x(t)$ is bounded, too. This completes the proof.
%Follow the analysis of the discrete-time network, we can get similar
%conclusion to the continuous-time network.

\begin{remark}
In \cite{Chen11'}, the sufficient condition to guarantee cluster consensus is that the number of clusters is equal to the period of agents.
The period of agent $i$ is the greatest common divisor of the lengths of paths starting form agent $i$ to itself.
To apply the results in \cite{Chen11'}, the period of all agents should be no less than 2.
In our paper, we assume the existence of self-links, which means the period of every agent is 1. So, the results in \cite{Chen11'} cannot be
employed in our situation.
\end{remark}

\section{Numerical Examples}

Cluster consensus is a new issue in the coordination
control. Despite that a huge number of papers were concerned with
complete consensus, there are a small amount of papers involved with
cluster consensus. Moreover, all of them cannot handle the scenario
in the paper. For example, \cite{Yu09} and \cite{Yu09'} investigated group
consensus in continuous-time network with fixed and switching
topologies respectively. Instead, in our paper, we study the
discrete-time network. Even though \cite{Chen11'} investigated the cluster
consensus in discrete-time network, it was concluded that cluster
consensus can be achieved if the graph topology is fixed and
strongly connected and the number of clusters equals to the period
of agents. Hence, the period of agents should be larger than 1. But
in our paper, the assumption that each agent has self-link means
that the period of agents in our algorithm is 1. For these reasons,
their results can hardly be applied to our case.

In this section, we provide an application example by a modified
non-Bayesian social learning model. Social learning can be described
as the process by which individuals infer information about some
alternative by observing the choices of others. In
\cite{Jadbabaie11}, a new social learning model was proposed, by
which an individual updates his/her belief as a convex combination of
the Bayesian posterior beliefs based on its private signal and the
beliefs of its neighbors at the previous time. In details, let
$\Theta=\{\theta_{1},\cdots, \theta_{m}\}$ denote a finite set of
possible states of the world and $\mu_{i,t}(\theta)$ denote the
probability (belief in their terminology) of individual $i$ about
state $\theta\in \Theta$ at time $t$. Conditional on the state
$\theta$, a signal vector
$\omega_{t}=(\omega_{1,t},\cdots,\omega_{n,t})\in S_{1}\times \cdots
\times S_{n}$ is generated by the likelihood function $l(\cdot
|\theta)$, where signal $\omega_{i,t}$ is the signal privately
observed by agent $i$ at period $t$ and $S_{i}$ denotes the signal
space of agent $i$. $l_{i}(\cdot |\theta)$ is the $i$-th
marginal of $l(\cdot |\theta)$. It is assumed that every agent $i$ knows this
conditional likelihood function. The one-step-ahead forecast of
agent $i$ at time $t$ is given by
$m_{i,t}(\omega_{i,t+1})=\sum_{\theta\in
\Theta}l_{i}(\omega_{i,t+1}|\theta)\mu_{i,t}(\theta)$.
The $k$-step-ahead forecast of agent $i$ at time $t$ is similarly given by
$m_{i,t}(\omega_{i,t+1},\cdots,\omega_{i,t+k})=\sum_{\theta\in
\Theta}(\prod_{r=1}^{k}l_{i}(\omega_{i,t+r}|\theta))\mu_{i,t}(\theta)$
Then, the
belief updating rule can be written as

\begin{eqnarray}
\mu_{i,t+1}(\theta)=a_{ii}\mu_{i,t}(\theta)
\frac{l_{i}(\omega_{i,t+1}|\theta)}{m_{i,t}(\omega_{i,t+1})}+\sum_{j\in
\mathcal{N}_{i}}a_{ij}\mu_{j,t}(\theta)\label{LA}
\end{eqnarray}

\cite{Jadbabaie11} considered the case that each agent may face an
identification problem in the sense that agent may not be able to
distinguish between two states. Observationally equivalence is used
to reflect the identification problem. Two states are
observationally equivalent from the point of view of agent $i$, if
the conditional distributions of agent $i$'s signals under the two
states coincide. As proved in \cite{Jadbabaie11}, all briefs
asymptotically coincide by this algorithm. This confirms the facts
that the interaction among individuals can eliminate the initial
difference among them and converge to an agreement.

For any state $\theta$, (\ref{LA}) can be rewritten in matrix form:
\begin{eqnarray}
\mu_{t+1}(\theta)=A\mu_{t}(\theta)+e_{t}(\theta)\label{Matrix form}
\end{eqnarray}
here $e_{t}(\theta)=(e_{1,t}(\theta),\cdots,e_{n,t}(\theta))^{\top}$ and $e_{i,t}(\theta)=a_{ii}\mu_{i,t}(\theta)
(\frac{l_{i}(\omega_{i,t+1}|\theta)}{m_{i,t}(\omega_{i,t+1})}-1)$.
For state $\hat{\theta}$ that is observationally equivalent
to $\theta^{*}$, the one-step-ahead forecasts and $k$-step-ahead forecasts
respectively satisfy
\begin{eqnarray*}
m_{i,t}(\omega_{i,t+1})\rightarrow l_{i}(\omega_{i,t+1}|\hat{\theta}),  ~~~t\rightarrow \infty
\end{eqnarray*}
and
\begin{eqnarray}\label{forecast}
m_{i,t}(\omega_{i,t+1},\cdots,\omega_{i,t+k})\rightarrow
\prod_{r=1}^{k}l_{i}(\omega_{i,t+r}|\hat{\theta}),~~~t\rightarrow \infty
\end{eqnarray}

Therefore, $e_{i,t}(\hat{\theta})$ converges to zero almost surely
as time goes on. Then from matrix and probability theories, the
existence of $\lim_{t\rightarrow \infty}\mu_{i,t}(\hat{\theta})$ can
be obtained. For state $\theta$ that is not observationally
equivalent to $\theta^{*}$, there exist a positive integer
$\hat{k}_{i}$, a sequence of signals
$(\hat{s}_{i,1},\cdots,\hat{s}_{i,\hat{k}_{i}})$ and constant
$\delta_{i}\in (0,1)$ such that $\prod_{r=1}^{\hat{k}_{i}}
\frac{l_{i}(\hat{s}_{i,r}|\theta)}{l_{i}(\hat{s}_{i,r}|\theta^{*})}\leq
\delta_{i}$, combining with the $k$-step-ahead forecast
(\ref{forecast}), $\mu_{i,t}(\theta)\rightarrow 0$ a.s. can be
obtained.

Here, we assume that all states $\theta_{j}\in\Theta$ are
observationally equivalent for all individuals. Under this
assumption,
$\frac{l_{i}(\omega_{i,t+1}|\theta_{j})}{m_{i,t}(\omega_{i,t+1})}=1$
always are true. This implies that the signals observed have no
effect in this situation,  thus we remove the conditional likelihood
term in (\ref{LA}). In addition, we consider that  the belief of
each individual is affected by different religious beliefs or
cultural backgrounds. This affection flags the sub-group that each
individual belong to. Consider the group with $9$ individuals that
are divided into three groups: $\mathcal{C}_{1}=\{1,2,3\}$,
$\mathcal{C}_{2}=\{4,5,6\}$ and $\mathcal{C}_{3}=\{7,8,9\}$. We
denote auxiliary terms, $I_{i}(t)$, as the external inputs to the
learning model (\ref{LA}), in order to denote the influence of the
religious beliefs and/or cultural backgrounds and they are different
with respect to sub-groups (cluster). These terms are regarded as
the flags that distinguish the different sub-groups (clusters).
Hence, the dynamic model (\ref{LA}) becomes:
\begin{eqnarray}\label{simulation}
\mu_{i,t+1}(\theta)=a_{ii}\mu_{i,t}(\theta) +\sum_{j\in
\mathcal{N}_{i}}a_{ij}\mu_{j,t}(\theta) +I_{i}(t)
\end{eqnarray}
with the cultural/religious terms:
\begin{eqnarray*}
I_{i}(t)=cu(t)\sigma_{k}(\theta),~i\in\mathcal C_{k},~k=1,2,3,
\end{eqnarray*}
where $c$ denotes the influence strength. To guarantee
$\mu_{i,t}(\theta)\in [0,1]$, we assume the inter-cluster nonidentical input $u(t)$ is
periodic with a period $T=2$ and $u_{k}+u_{k+1}=0$. For every $i$
and $t$, to guarantee
$\mu_{i,t}(\theta_{1})+\mu_{i,t}(\theta_{2})=1$,  we demand
$\sigma_{i}(\theta_{1})+\sigma_{i}(\theta_{2})=0$. It can be seen
that the modified social learning  model (\ref{simulation}) is a
special case of the model (\ref{Mu-uniformCondition}).

To illustrate the availability of our results, we consider the state
space has two states: $\Theta=\{\theta_{1},\theta_{2}\}$. The
coupling matrix $A=[a_{ij}]$ satisfies the inter-cluster influence
condition, and suppose $\{k|\mathcal{N}_{i}\cap\mathcal{C}_{k}\neq
\emptyset\}$ is identical to all $i\in \mathcal{C}_{p},p=1,2,3$.
Denote $d_{iq}$ the number of agents in set
$\mathcal{N}_{i}\cap\mathcal{C}_{q}$ and for $q\in
\{k|\mathcal{N}_{i}\cap \mathcal{C}_{k}\neq \emptyset\}, j\in
\mathcal{N}_{i}\cap \mathcal{C}_{q}$, take
$a_{ij}=\frac{\beta_{pq}}{d_{iq}}$. For any $p$ and any $q \in
\{k|\mathcal{N}_{i}\cap\mathcal{C}_{k}\neq \emptyset\}$, $\sum_{j\in
\mathcal{C}_{q}}\frac{\beta_{p,q}}{d_{iq}}=\sum_{j\in
\mathcal{C}_{q}}\frac{\beta_{p,q}}{d_{i'q}}=\beta_{p,q}$ always
holds for $\forall i,i'\in \mathcal{C}_{p}$, i.e. the coupling
matrix in (\ref{simulation}) has the common inter-cluster influence.
We use $B=[\beta_{pq}]_{p,q=1}^{3}$ to inflect the inter-cluster
influence among clusters, and choose $u(2l)=-u(2l+1)=1$, for all
$l\in \mathbb{N}$.

\begin{figure*}[!t]
\begin{center}
\includegraphics[width=.8\textwidth]{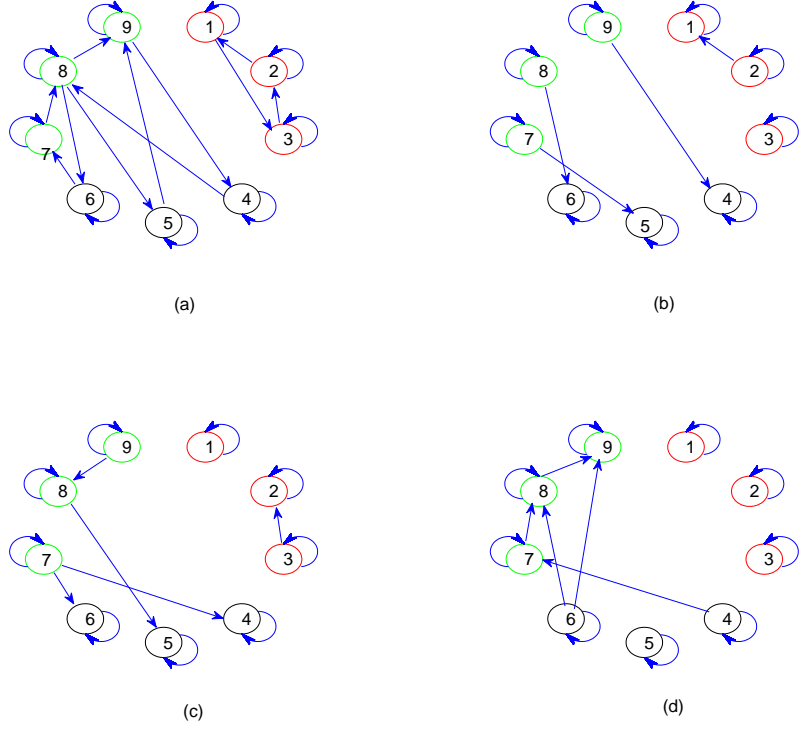}
\end{center}
\caption{All of these graphs have self -links. Example A simulate
the network with fixed topology (a) and example B simulate the
network with topologies switching in (b),(c),(d)}.\label{graph}
\end{figure*}

%We use
%\begin{eqnarray*}
%\eta(\mu_{t}(\theta_{j}))=\max_{i\in\mathcal{C}_{p},j\in\mathcal{C}_{q},p\neq q}|x_{i}(t)-x_{j}(t)|
%\end{eqnarray*}
%to measure the difference between clusters and
%\begin{eqnarray*}
%\Delta_{\mathcal{C}}(x(t))=\max_{p}\max_{i,i'\in\mathcal{C}_{p}}|x_{i}(t)-x_{i'}(t)|
%\end{eqnarray*}
%to measure the difference of states in the same clusters.

\textbf{A. Static topology}

In this example, the graph is depicted in Fig \ref{graph} (a). We
take the matrix $B$ as:
\begin{eqnarray*}
B=\left[
\begin{array}{ccc}
  1 & 0 & 0 \\
  0& 1/2 & 1/2 \\
  0 & 1/2 & 1/2
\end{array}
\right]
\end{eqnarray*}
and can see that the graph has cluster spanning trees and the roots
of groups $\mathcal{C}_{1,2,3}$ are $3$, $7$ and $7$ respectively.
Therefore, all conditions in Theorem 1 hold. Then (\ref{simulation})
reaches cluster consensus generically.
\begin{figure}[!t]
\begin{center}
\includegraphics[width=.5\textwidth]{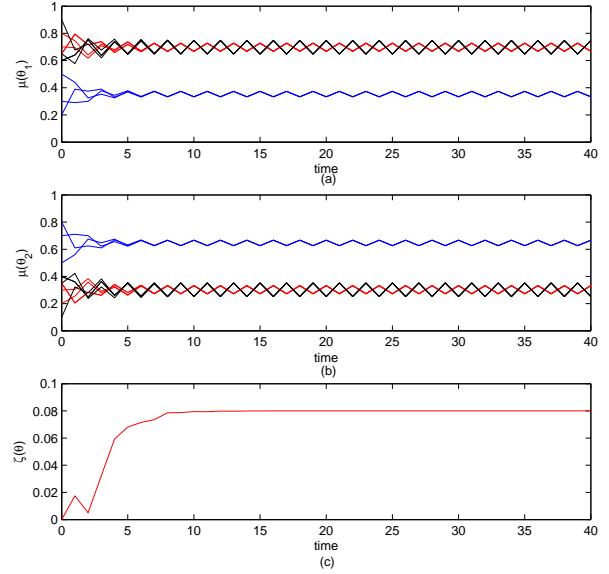}
\end{center}
\caption{The dynamical behavior of beliefs $\mu(\theta_{1})$,
$\mu(\theta_{2})$ and $\zeta(\theta)$ with respect to example A with
randomly picked initial values. In (a) and (b), blue, red and black
curves show the dynamical behaviors of individuals in group
$\mathcal{C}_{1},\mathcal{C}_{2},\mathcal{C}_{3}$ respectively.}
\label{se}
\end{figure}
The dynamical behaviors of the beliefs $\mu_{i,t}(\theta_{j}),
i=1,\cdots,9, j=1,2$ are shown in Fig \ref{se} (a) and (b). It is
clear that they are asymptotically convergent, which means different
groups of individuals can realize intra-cluster synchronization. In
Fig \ref{se} (c), the dynamical behaviors of
$\zeta(\theta_{j})=|\mu_{\mathcal{C}_{2}}(\theta_{j})
-\mu_{\mathcal{C}_{3}}(\theta_{j})|, j=1$ is plotted and it does not
converge to zero, which means that although groups $\mathcal C_{2}$
and $\mathcal C_{3}$ are strongly connected, the influence of
different religious beliefs or cultural backgrounds still cannot be
ignored.

\textbf{B. Switching topologies}\\
In this example, the graph topology is switching among the
topologies given in Fig \ref{graph} (b), (c) and (d) periodically.
Noting that none of these graphs has cluster spanning trees, i.e.
the condition in Theorem \ref{intra-cluster} does not hold. However,
the union graph of those in Fig \ref{graph} (b), (c) and (d) has
cluster spanning trees and the roots of groups $\mathcal{C}_{1,2,3}$
are agents $3$, $7$ and $7$ respectively. We pick an identical
matrix $B$ w.r.t. the clustering for the three graphs as
\begin{eqnarray*}
B=\left[
\begin{array}{ccc}
  1 & 0 & 0 \\
  0& 1/2 & 1/2 \\
  0 & 1/2 & 1/2
\end{array}
\right].\end{eqnarray*}

\begin{figure}[!t]
\begin{center}
\includegraphics[width=.5\textwidth]{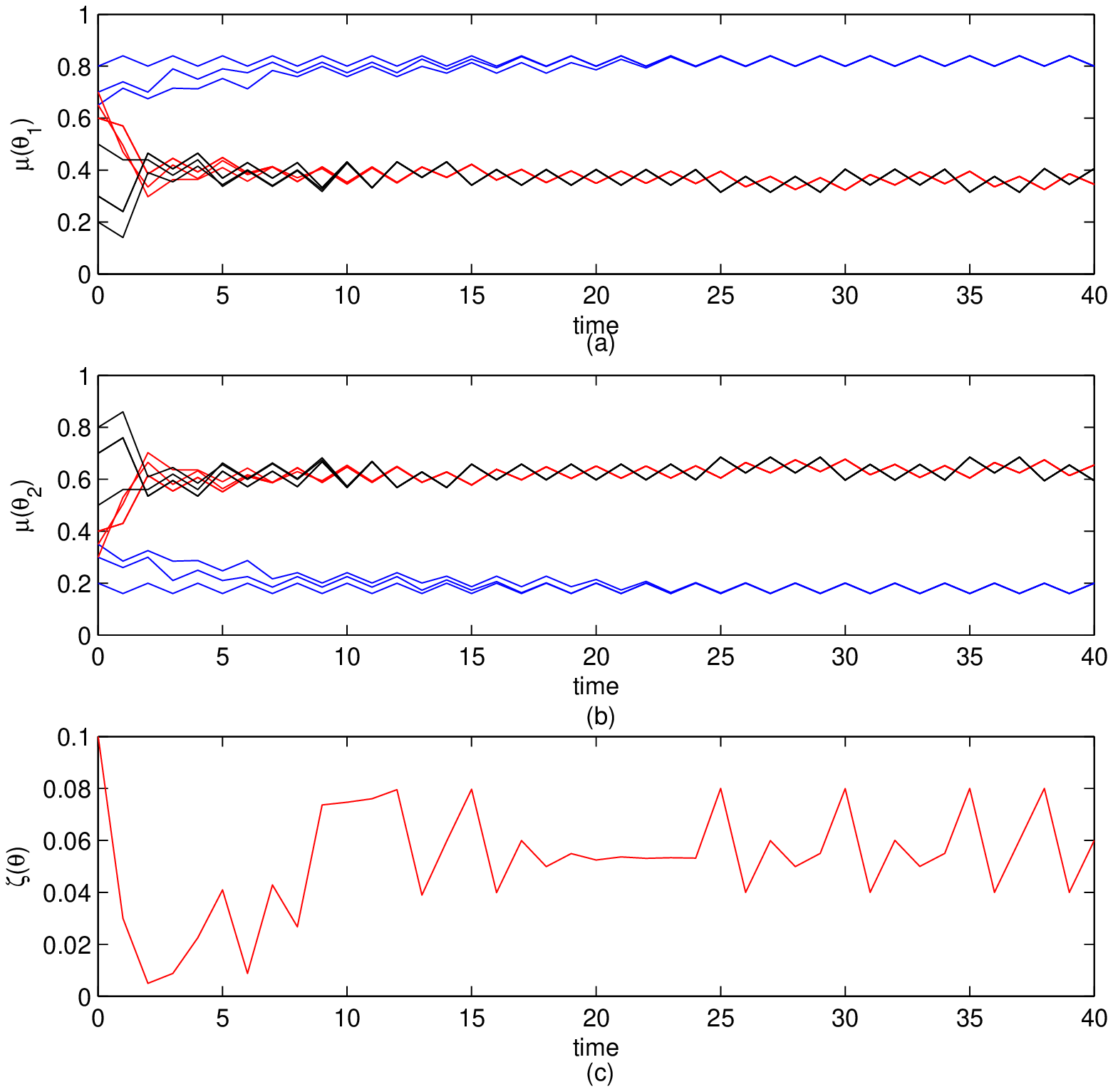}
\end{center}
\caption{The dynamical behavior of states $\mu(\theta_{1})$,
$\mu(\theta_{2})$ and $\zeta(\theta)$ with respect to example B with
randomly picked initial values. In (a) and (b), blue, red and black
curves show the dynamical behaviors of individuals in group
$\mathcal{C}_{1},\mathcal{C}_{2},\mathcal{C}_{3}$ respectively.}
\label{xe}
\end{figure}

Hence, all assumptions in Theorem 4 hold. Therefore,
(\ref{simulation}) with switching topologies can achieve cluster
consensus. The dynamical behaviors of beliefs
$\mu_{i,t}(\theta_{j}), 1 \leq i\leq 9$ are shown in Fig \ref{xe}
(a) and (b), the dynamics of
$\zeta(\theta_{j})=|\mu_{\mathcal{C}_{2}}(\theta_{j})
-\mu_{\mathcal{C}_{3}}(\theta_{j})|, j=1$ is plotted in Fig \ref{xe}(c)
respectively. All of them show that the cluster consensus is
perfectly reached and $\mu_{i,t}(\theta_{j}), 1 \leq i\leq 9$ is
convergent.

\begin{figure}[!t]
\begin{center}
\includegraphics[width=.5\textwidth]{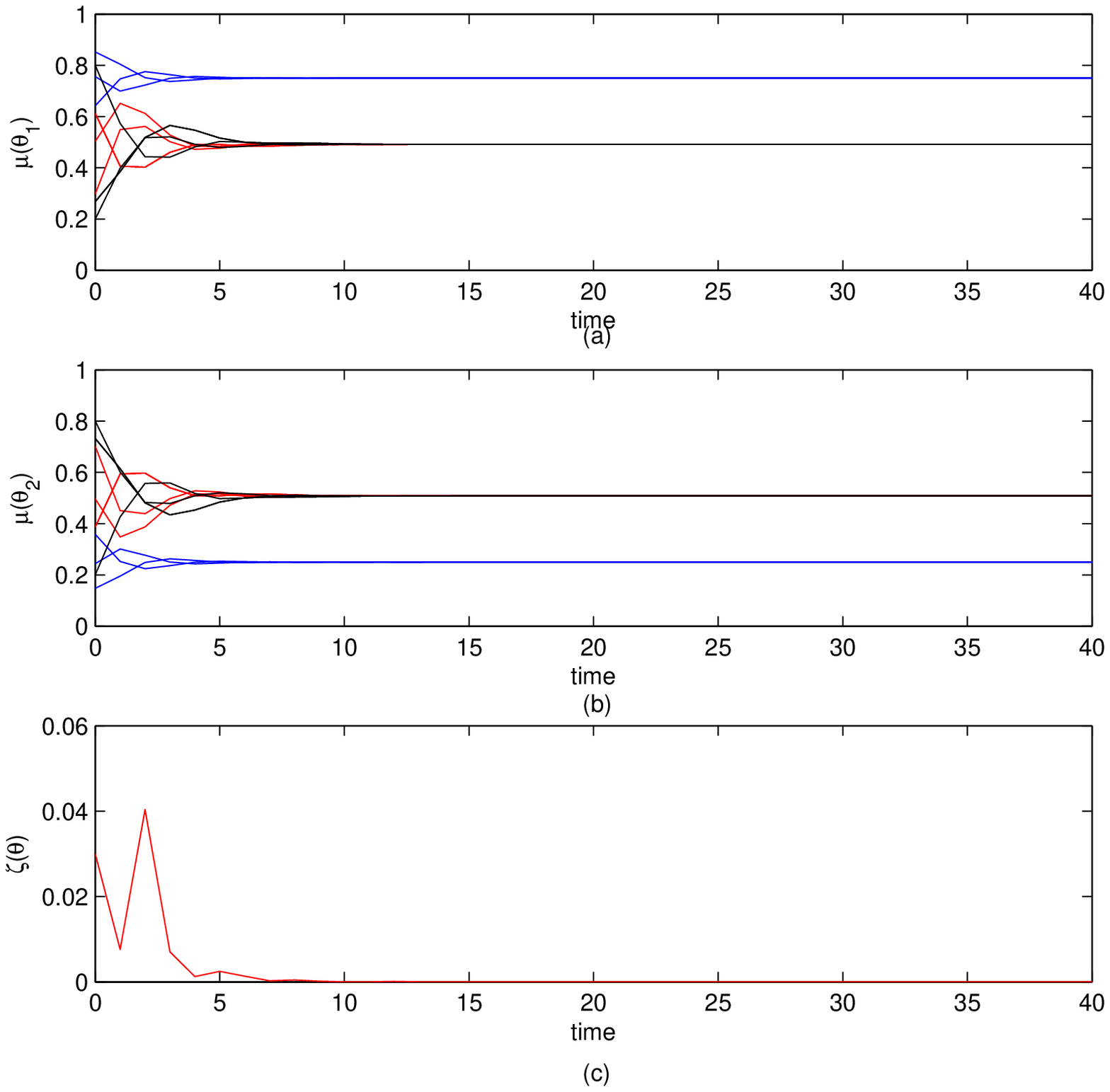}
\end{center}
\caption{The dynamical behavior of states $\mu(\theta_{1})$,
$\mu(\theta_{2})$ and $\zeta(\theta)$ with respect to example A with
randomly picked initial values. In (a) and (b), blue, red and black
curves show the dynamical behaviors of individuals in group
$\mathcal{C}_{1},\mathcal{C}_{2},\mathcal{C}_{3}$ respectively.}
\label{c1}
\end{figure}

Now, to better illustrate the role of the inter-cluster nonidentical inputs, we give a
simulation based on (\ref{simulation}) without inputs, see Fig
\ref{c1}. The dynamical behaviors of beliefs $\mu_{i,t}(\theta_{j}),
i=1,\cdots,9, j=1,2$ are shown in Fig \ref{c1} (a) and (b). In  Fig
\ref{c1} (c), the dynamical behavior of
$\zeta(\theta_{j})=|\mu_{\mathcal{C}_{2}}(\theta_{j})
-\mu_{\mathcal{C}_{3}}(\theta_{j})|$ is plotted, which means the
groups $\mathcal C_{2}$ and $\mathcal C_{2}$ cannot separate.
Compare with Fig \ref{se}(c), we can see that the inter-cluster nonidentical inputs
play key roles in separating different groups.

\section{Conclusions}

The idea for studying consensus of multi-agent systems sheds light
on cluster consensus analysis. In this paper, we study cluster
consensus of multi-agent systems via inter-cluster nonidentical inputs. We derive the
criteria for cluster consensus in both discrete-time systems with
fixed or switching graph topologies. The difference between
clustered states are guaranteed by the different inputs to
different clusters. We present if every cluster in the graph
corresponding to the system has a spanning tree, then the
multi-agent system reaches cluster consensus. The analysis is
presented rigorously based on algebraic graph theory and matrix
theory.  We use a modified non-Bayesian social learning model to
illustrate our theoretical results. In this model, the briefs of
individuals are described as the probability for the states and
updated by an interacted algorithms. We add an auxiliary term to
flag the difference of culture and/or region of different group of
individuals. The numerical results show that the social learning
algorithm can guarantee that the briefs of individuals in the same
cluster converge but the difference between any pair of groups,
owing to the auxiliary external input terms, permanently exists that
cannot be eliminated by the interactions.

\section*{Appendix}
\textbf{Proof of Lemma \ref{scrambling}:} For each cluster $\mathcal
C_{p}$ and each pair of vertices $v_{p_{1}},~v_{p_{2}}\in\mathcal
C_{p}$, let $\mathcal V_{1}^{t}$ and $\mathcal V_{2}^{t}$ be the
neighborhoods to $v_{p_{1}}$ and $v_{p_{2}}$ respectively in the
graph $\mathcal G(A_{1}^{t})$. The fact that each $A(t)$ has all
nodes self-linked implies that $\mathcal V_{i}^{t}\subset\mathcal
V_{i}^{t+1}, i=1,2$ respectively. In the following, we are going to
prove that $\mathcal
V_{1}^{t}\bigcap\mathcal V_{2}^{t}\ne\emptyset$ holds for at least some $t\le n$.

If $t<n,$ $\mathcal V_{1}^{t}\bigcap\mathcal V_{2}^{t}=\emptyset$,
then $\#[\mathcal V_{1}^{t}\bigcup\mathcal V_{2}^{t}]\ge t+1$.

We will prove it by induction. By
the assumptions, there is a cluster root in $\mathcal G (A(1))$ that
has paths towards the vertices $v_{p_{1}}$ and $v_{p_{2}}$, both
$\mathcal V_{1}^{1}$ and $\mathcal V_{2}^{1}$ are not empty. If
$\mathcal V_{1}^{1}\bigcap\mathcal V_{2}^{1}=\emptyset$, then
$\#[\mathcal V_{1}^{1}\bigcup\mathcal V_{2}^{1}]\ge 2$.

Suppose $\mathcal V_{1}^{t}\bigcap\mathcal V_{2}^{t}=\emptyset$ and
$\#[\mathcal V_{1}^{t}\bigcup\mathcal V_{2}^{t}]\ge t+1$. We will
prove $\#\bigg(\mathcal V_{1}^{t+1}\bigcup\mathcal
V_{2}^{t+1}\bigg)\ge t+2$.

In fact, let $v_{1}$ be the root vertex in the graph $\mathcal
G(A(t+1))$  having paths towards $v_{p_{1}}$ and $v_{p_{2}}$. We
select their shortest paths:
$(v_{k_{1}},v_{k_{2}},\cdots,v_{k_{P}})$ and
$(v_{l_{1}},v_{l_{2}},\cdots,v_{l_{Q}})$, from $v_{1}$ to
$v_{p_{1}}$ and $v_{p_{2}}$ respectively, with
$v_{k_{1}}=v_{l_{1}}=v_{1}$, $v_{k_{P}}=v_{p_{1}}$ and
$v_{l_{Q}}=v_{p_{2}}$. If one of the paths has one vertex not
belonging to the corresponding $\mathcal V_{1}^{t}$ or $\mathcal
V_{2}^{t}$. Without loss of generality, we assume that
$(v_{k_{1}},v_{k_{2}},\cdots,v_{k_{P}})$ has vertices not belonging
to $\mathcal V_{1}^{t}$ and let $v_{k_{r_{0}}}$ be the index such
that
\begin{itemize}
\item for each $r>r_{0}$, $v_{k_{r}}\in\mathcal V_{1}^{t}$;
\item $v_{k_{r_{0}}}\notin\mathcal V_{1}^{t}$.
\end{itemize}
This implies that
\begin{eqnarray*}
[A_{1}^{t+1}]_{v_{k_{r_{0}}},v_{k_{P}}}\ge [A(t+1)]_{k_{r_{0}},k_{r_{0}+1}}[A_{1}^{t}]_{v_{k_{r_{0}+1}},v_{k_{P}}}>0
\end{eqnarray*}
holds. This implies that $v_{k_{r_{0}}}\in\mathcal V_{1}^{t+1}$.
Hence,
\begin{eqnarray*}
\#\bigg(\mathcal V_{1}^{t+1}\bigcup\mathcal
V_{2}^{t+1}\bigg)\ge\#\bigg(\mathcal V_{1}^{t}\bigcup\mathcal
V_{2}^{t}\bigg)+1\ge t+2.
\end{eqnarray*}
Thus, either for some $t<n$, $\mathcal V_{1}^{t}\bigcap\mathcal
V_{2}^{t}\ne\emptyset$ or
\begin{eqnarray*}
\#\bigg(\mathcal V_{1}^{n}\bigcup\mathcal
V_{2}^{n}\bigg)\ge n+1
\end{eqnarray*}
which implies $V_{1}^{n}\bigcap\mathcal V_{2}^{n}\ne\emptyset$.
Therefore, there exists some $t\le n$ such that $\mathcal
V_{1}^{t}\bigcap\mathcal V_{2}^{t}\ne\emptyset$. Proof of the lemma
is completed.

\textbf{Proof of Claim 1:}
\begin{eqnarray*}
\mathbb R^{n}=\mathcal S_{\mathcal C}+V.
\end{eqnarray*}

For this purpose, we define a nonsingular matrix
$P=[P_{1},\cdots,P_{n}]\in\mathbb R^{n,n}$ such that the first $K$
column vectors compose a basis of $\mathcal S_{\mathcal C}$. In
particular, we chose each $P_{k}$, $k=\onetoK$, as
\begin{eqnarray*}
[P_{k}]_{i}=\begin{cases}1&i\in\mathcal C_{k}\\
0&otherwise.\end{cases}
\end{eqnarray*}
Define
\begin{eqnarray*}
\hat{A}(t)\triangleq P^{-1}A(t)P=\left[\begin{array}{ll}\hat{A}_{1,1}&\hat{A}_{1,2}(t)\\
0&\hat{A}_{2,2}(t)\end{array}\right],
\end{eqnarray*}
where the bottom-left block equals to zero since the subspace $\mathcal S_{\mathcal C}$ is invariant by $A(t)$ and the top-left block $\hat{A}_{1,1}$ is a static matrix due to $\mathcal B_{3}^{*}$. Furthermore, since all eigenvalues of $B$, defined in (\ref{B1}), of which the modules equal to $1$  should equal to $1$, owing to the fact that all matrices $A(t)$ have all diagonal elements positive, we can select $Q_{1}$ with the first several columns composing of the basis of the eigenspace of the static sub-matrix $\hat{A}_{1,1}$ corresponding to eigenvalue $1$ and all last $n-K$ columns was chosen to guarantee $Q_{1}$ is nonsingular. Construct a new linear transformation  $Q$ has the form as:
\begin{eqnarray*}
Q=\left[\begin{array}{cc}Q_{1}&0\\
0&I_{n-K}\end{array}\right].
\end{eqnarray*}
Then, we further make linear transformation with $Q$ over $\hat{A}(t)$ resulting in:
\begin{eqnarray*}
\tilde{A}(t)\triangleq Q^{-1}\hat{A}(t)Q=\left[\begin{array}{ll}\tilde{A}_{1,1}&\tilde{A}_{1,2}(t)\\
0&\tilde{A}_{2,2}(t)\end{array}\right],
\end{eqnarray*}
where $\tilde{A}_{1,1}$ has the following block form:
\begin{eqnarray*}
\tilde{A}_{1,1}=\left[\begin{array}{cc}\tilde{A}_{1,1}^{1,1}&0\\
0&\tilde{A}_{1,1}^{2,2}\end{array}\right].
\end{eqnarray*}
with all eigenvalues of $\tilde{A}_{1,1}^{1,1}$ equal to $1$ and $\rho(\tilde{A}_{1,1}^{2,2})<1$. Accordingly, we write
\begin{eqnarray*}
\tilde{A}_{1,2}(t)=\left[\begin{array}{c}
\tilde{A}_{1,2}^{1}(t)\\
\tilde{A}_{1,2}^{2}(t)\end{array}\right].
\end{eqnarray*}

Thus, we define
\begin{eqnarray*}
\tilde{A}_{0}^{t}=\left[\begin{array}{ll}(\tilde{A}_{1,1})^{t+1}&\tilde{A}_{1,2}^{(t)}\\
0&(\tilde{A}_{2,2})_{0}^{t}\end{array}\right]
\end{eqnarray*}
where $(\cdot)_{0}^{t}$ denotes the left matrix product from $0$ to $t$, as defined before.

We define the {\em projection radius} (w.r.t. $\mathcal C$) of $A(t)$ as follows:
\begin{eqnarray*}
\rho_{\mathcal C}(A(\cdot))=\overline{\lim}_{t\to\infty}\bigg\{\|(\tilde{A}_{2,2})_{0}^{t-1}\|\bigg\}^{1/t}
\end{eqnarray*}
and the {\em cluster Hajnal diameter} (w.r.t. $\mathcal C$) of $A(t)$ as follows:
\begin{eqnarray*}
\Delta_{\mathcal C}(A(\cdot))=\overline{\lim}_{t\to\infty}\bigg\{\|\Delta_{\mathcal C}(A_{0}^{t-1})\|\bigg\}^{1/t}
\end{eqnarray*}
for some norm $\|\cdot\|$ that is induced by vector norm. It can be seen that the projection radius and cluster Hajnal diameter are independent of the selection of the matrix norm and the matrix $P$.
First, we shall prove that the projection radius equals to the Hajnal diameter.
\begin{lemma}
$\rho_{\mathcal C}(A(\cdot))=\Delta_{\mathcal C}(A(\cdot))$.
\end{lemma}
{\em Proof.} The proof is quite similar to that in \cite{Lu..} and can be regarded as a generalization of Lemma 2.4 in \cite{Lu..}. For any $d>\rho_{\mathcal C}(A(\cdot))$, there exists $T>0$ such that the inequality
\begin{eqnarray*}
\|(\tilde{A}_{2,2})_{0}^{t-1}\|<d^{t}
\end{eqnarray*}
for all $t>T$. Then
\begin{eqnarray*}
&&\bigg\|\tilde{A}_{0}^{t-1}-\left[\begin{array}{c}E_{K}\\0\end{array}\right]\left[(\tilde{A}_{1,1})^{t-1},\tilde{A}_{1,2}^{(t-1)}\right]\bigg\|\\
&&=\bigg\|\left[\begin{array}{ll}0&0\\
0&(\tilde{A}_{2,2})_{0}^{t-1}\end{array}\right]\bigg\|\le C d^{t}
\end{eqnarray*}
for some $C>0$ and all $t>T$. Thus,
\begin{eqnarray*}
\bigg\|A_{0}^{t-1}-P\left[\begin{array}{c}E_{K}\\0\end{array}\right]\left[(\tilde{A}_{1,1})^{t-1},\tilde{A}_{1,2}^{(t-1)}\right]P^{-1}\bigg\|\le C_{1}d^{t},
\end{eqnarray*}
for some $C_{1}>0$ and all $t>T$. Let
\begin{eqnarray*}
&&G=P\left[\begin{array}{c}E_{K}\\0\end{array}\right]=[P_{1},\cdots,P_{K}],\\
&&H=\left[(\tilde{A}_{1,1})^{t-1},\tilde{A}_{1,2}^{(t-1)}\right]P^{-1}.
\end{eqnarray*}
Since each $P_{k}\in \mathcal S_{\mathcal C}$ for all $k=\onetoK$, each column vector in the matrix $G\cdot H$ should belong to $\mathcal S_{\mathcal C}$, too.  So, according to the definition of Hajnal diameter, we have
\begin{eqnarray*}
\Delta_{\mathcal C}(A_{0}^{t-1})\le 2 C_{1} d^{t}
\end{eqnarray*}
for all $t\ge T$. This implies that $\Delta_{\mathcal C}(A(\cdot))\le d$. According to the arbitrariness of $d$, we have $\Delta_{\mathcal C}(A(\cdot))\le\rho_{\mathcal C}(A(\cdot))$.

On the other hand, for any $d>\Delta_{\mathcal C}(A(\cdot))$, there exists $T>0$ such that $\Delta_{\mathcal C}(A_{0}^{t-1})<d^{t}$ holds for all $t>T$. Without loss of generality, we suppose that the clustering $\mathcal C$ is successive, i.e., $\mathcal C_{1}=\{1,2,\cdots,n_{1}\}$, $\mathcal C_{2}=\{n_{1}+1,n_{1}+2,\cdots,n_{2}\}$,$\cdots$, $\mathcal C_{K}=\{n_{k-1}+1,n_{k-1}+2,\cdots,n_{K}\}$ with $n_{K}=n$. Select one single row in $A_{0}^{t-1}$ from each cluster and compose them into a matrix, denoted by $H$. Let $G=[P_{1},\cdots, P_{K}]$. Then the rows of $G\cdot H$ corresponding to the same cluster is identical. So,
\begin{eqnarray*}
\|A_{0}^{t-1}-G\cdot H\|\le C_{2}d^{t}
\end{eqnarray*}
for some $C_{2}>0$ and $t>T$. Then,
\begin{eqnarray*}
\|P^{-1}A_{0}^{t-1}P-P^{-1}G\cdot HP\|\le C_{3}d^{t}
\end{eqnarray*}
i.e.,
\begin{eqnarray*}
\|\left[\begin{array}{ll}(\hat{A}_{1,1})^{t-1}&\hat{A}_{1,2}^{(t-1)}\\
0&(\hat{A}_{2,2})_{0}^{t-1}\end{array}\right]-\left[\begin{array}{ll}Y&Z\\
0&0\end{array}\right]\|\le C_{3}d^{t}
\end{eqnarray*}
for some matrices $Y$ and $Z$, $C_{3}>0$ and all $t>T$. This implies that $\|(\hat{A}_{2,2})_{0}^{t-1}\|\le C_{4}d^{t}$ holds for some $C_{4}>0$ and all $t>T$. It can be seen that $(\hat{A}_{2,2})_{0}^{t-1}=(\tilde{A}_{2,2})_{0}^{t-1}$. Therefore, $\rho_{\mathcal C}(A(\cdot))\le d$. The arbitrariness of $d$ can guarantee $\Delta_{\mathcal C}(A(\cdot))\ge\rho_{\mathcal C}(A(\cdot))$. From both sides, we have $\Delta_{\mathcal C}(A(\cdot))=\rho_{\mathcal C}(A(\cdot))$.  This completes the proof of this lemma.

From Theorem  \ref{tv intra-cluster}, we can conclude $\Delta_{\mathcal C}(A(\cdot))<1$. Thus, $\rho_{\mathcal C}(A(\cdot))<1$.
For any $n$-dimensional vector $w_{0}$, we can write it as:
\begin{eqnarray*}
w_{0}=\left[\begin{array}{c}z_{0}\\u_{0}\\v_{0}\end{array}\right]
\end{eqnarray*}
where $z_{0}$ corresponds to the sub-matrix $\tilde{A}_{1,1}^{1,1}$,  $u_{0}$ corresponds to the sub-matrix $\tilde{A}_{1,1}^{2,2}$ and $v_{0}\in\mathbb R^{n-K}$. We rewrite $w_{0}$ as a sum of $w_{0}^{1}+w_{0}^{2}$ with
\begin{eqnarray*}
w_{0}^{1}=\left[\begin{array}{c}z_{0}^{1}\\0\\0\end{array}\right],~w_{0}^{2}=\left[\begin{array}{c}
z_{0}^{2}\\u_{0}\\v_{0}\end{array}\right]
\end{eqnarray*}
where $z_{0}^{1}+z_{0}^{2}=z_{0}$ that will be determined in the following. It is clear that $PQ w_{0}^{1}$ corresponds a vector in $\mathcal S_{C}$. So, if we could pick a suitable $z_{0}^{2}$ such that $\lim_{t\to\infty}(\tilde{A})_{0}^{t}w_{0}^{2}=0$, that is, $PQ w_{0}^{2}$ corresponds a vector in $V$. Therefore, for any $n$-dimensional vector $x_{0}$, we can find some $w_{0}$, such that $x_{0}=PQ w_{0}=PQw_{0}^{1}+PQw_{0}^{2}\in \mathcal{S}_{\mathcal{C}}+V$. This could complete the proof of the claim.

For this purpose, we consider the  following linear system:
\begin{eqnarray*}
\tilde{w}(t+1)=\tilde{A}(t)\tilde{w}(t),~\tilde{w}(0)=w_{0}^{2},
\end{eqnarray*}
which can be rewritten as the following component-wise form:
\begin{eqnarray*}
&&\begin{cases}
\tilde{w}_{1}(t+1)=\tilde{A}_{1,1}^{1,1}\tilde{w}_{1}(t)+\tilde{A}_{1,2}^{1}(t)\tilde{w}_{3}(t)\\
\tilde{w}_{2}(t+1)=\tilde{A}_{1,1}^{2,2}\tilde{w}_{2}(t)+\tilde{A}_{1,2}^{2}(t)\tilde{w}_{3}(t)\\
\tilde{w}_{3}(t+1)=\tilde{A}_{2,2}(t)\tilde{w}_{3}(t)
\end{cases}\\
&&{\rm with~}\tilde{w}_{1}(0)=z_{0}^{2},~\tilde{w}_{2}(0)=u_{0},~\tilde{w}_{3}(0)=v_{0}.
\end{eqnarray*}
It can be seen that $\lim_{t\to\infty}\tilde{w}_{3}(t)=0$ exponentially
because of $\rho_{\mathcal C}(A(\cdot))<1$ and $\lim_{t\to\infty}\tilde{w}_{2}(t)=0$ exponentially because of $\rho(\tilde{A}_{1,1}^{2,2})<1$ and the boundedness of $\tilde{A}_{1,2}^{2}(t)$. Without loss of generality, since $\rho_{\mathcal C}(A)<1$ and all eigenvalues of $(\tilde{A}_{1,1}^{1,1})^{-1}$ equal to $1$, for any $\epsilon_{0}\in(0,|\lambda_{2}|/2)$, we have $\|(\tilde{A}_{2,2})_{0}^{t}\|\le M_{2}\exp[-(|\lambda_{2}|-\epsilon_{0})t]$, $\|(\tilde{A}_{1,1}^{1,1})^{-1}\|<\exp(\epsilon_{0})$ and $\|\tilde{A}_{1,2}^{1}(t)\|\le M_{0}$ for some $M_{0}>0$, $\lambda_{0}>0$, all $t\ge 0$ and some norm $\|\cdot\|$. Note that
\begin{eqnarray*}
\tilde{w}_{1}(t)=(\tilde{A}_{1,1}^{1,1})^{t}z_{0}^{2}+\sum_{k=0}^{t}(\tilde{A}_{1,1}^{1,1})^{t-k}
\tilde{A}_{1,2}^{1}(k)(\tilde{A}_{2,2})_{0}^{k}v_{0}.
\end{eqnarray*}
Since
\begin{eqnarray*}
&&\|(\tilde{A}_{1,1}^{1,1})^{-k}\|\cdot\|\tilde{A}_{1,2}^{1}(k)\|\cdot\|(\tilde{A}_{2,2})_{0}^{k}\|\\
&&\le
\exp(\epsilon_{0} k-[|\lambda_{2}|-\epsilon_{0}] k)M_{2}^{2}\\
&&\le\exp(-[|\lambda_{2}|+2\epsilon_{0}]k)M_{2}^{2},
\end{eqnarray*}
we let
\begin{eqnarray*}
R=\sum_{k=0}^{\infty}(\tilde{A}_{1,1}^{1,1})^{-k}
\tilde{A}_{1,2}^{1}(k)(\tilde{A}_{2,2})_{0}^{k}
\end{eqnarray*}
of which the limit exists in the norm sense and the operator $R$ is well-defined. Let us consider a subspace of $\mathbb R^{n}$:
\begin{eqnarray*}
\tilde{V}=\bigg\{[z^{\top},u^{\top},v^{\top}]^{\top}\in\mathbb R^{n}:~z=-Rv\bigg\}.
\end{eqnarray*}
If we pick $z_{0}^{2}$ such that $w_{0}^{2}\in\tilde{V}$, then we have
\begin{eqnarray*}
&&(\tilde{A}_{1,1}^{1,1})^{-t}\tilde{w}_{1}(t)=z_{0}^{2}+\sum_{k=0}^{t}(\tilde{A}_{1,1}^{1,1})^{-k}
\tilde{A}_{1,2}^{1}(k)(\tilde{A}_{2,2})_{0}^{k}v_{0}\\
&&\rightarrow z_{0}^{2}+Rv_{0}=0
\end{eqnarray*}
exponentially as $t\to\infty$. So, $(\tilde{A})_{0}^{t}w_{0}^{2}$ converges to zero exponentially. This completes the proof.


\begin{thebibliography}{99}
\bibitem{Vidal03}
R. Vidal, O. Shakernia, and S. Sastry, Formation control of nonholonomic mobile robots omnidirectional visual servoing and motion segmentation.
{\it Proc. IEEE Conf. Robotics and Automation, 2003, pp. 584-589.}
\bibitem{Cortes03}
J. Cortes and F. Bullo, Coordination and geometric optimazation via distributed dynamical systems.
{\it SIAM J. Control Optim., May 2003.}
\bibitem{Fax04}
A. Fax and R. M. Murray, Information flow and cooperative control of vehicle formations.
{\it IEEE Trans. Automat. Contr., 49 (2004), pp. 1465-1476.}
\bibitem{Reynolds87}
C. W. Reynolds, Flocks, herds, and schools: A distributed behavioral model.
{\it Proc. Comp. Graphics ACM SIGGRAPH'87 Conf., 1987, vol. 21, pp. 25-34.}
\bibitem{Vicsek95}
T. Vicsek, A. Czirook, E. Ben-Jacob, I. Cohen, and O. Shochet, Novel type of phase transition in a system of self-derived particles.
{\it Phys. Rev. Lett., vol. 75, no. 6, pp. 1226-1229, 1995.}
\bibitem{Xiao04}
L. Xiao and S.Boyd, Fast linear iterations for distributed averaging.
{\it Syst. Control Lett., vol. 53, pp. 65-78, 2004.}
\bibitem{Acemoglu}
D. Acemoglu, M. A. Dahleh, I. Lobel, and A. Ozdaglar, Bayesian learning  in social  networks,
{\it Review of Economic Studies, no.1, pp. 1�C34, 2010.}

\bibitem{Jadbabaie11}
A. Jadbabaie, P. Molavi, A. Sandroni, and A. Tahbaz-Salehi, Non-
bayesian  social  learning.
{\it PIER Working Paper No.11-025, August 2011.}

\bibitem{Liu qipeng}
Q. P. Liu, A. L. Fang, L. Wang, and X. F. Wang, Non-Bayesian learning in social networks with time-varying weights,
{\it Proceedings of the 30th Chinese Control Coference, July 22-24, 2011.}

\bibitem{Ren05}
W. Ren and R. W. Beard, Consensus seeking in multiagent systems under dynamically changing interaction topologies.
{\it IEEE Trans. Automat. Control, 50(2005), pp. 655-661.}
\bibitem{Moreau05}
L. Moreau, Stability of multiagent systems with time-dependent communication links.
{\it IEEE Trans. Automat. Control, 50(2005), pp. 169-182.}
\bibitem{Porfiri07}
M. Porfiri and D. J. Stilwell, Consensus seeking over random weighted directed graphs.
{\it IEEE Trans. Automat. Control, 52(2007), pp. 1767-1773.}
\bibitem{Hatano05}
Y. Hatano and M. Mesbahi, Agreement over random networks.
{\it IEEE Trans. Automat. Control, 50(2005), pp. 1867-1872.}
\bibitem{Liu08}
B. Liu and T. P. Chen, Consensus in networks of multi-agents with cooperation and competition via stochastically switching topologies.
{\it IEEE Trans. Neural Netw., 19(2008), pp. 1967-1973.}

\bibitem{Olf07}
R. Olfati-Saber, J. A. Fax and R. M. Murray, Consensus and
cooperation in networked multi-agent systems,
{\it Proceedings of the IEEE, 95(2007), pp. 215-233}.

\bibitem{Ren04}
W. Ren and R. W. Beard, Consensus of information under dynamically changing interaction topologies.
{\it Amer. Control Conf., 2004, pp. 4939-4944.}

\bibitem{Lu09}
W. Lu, F. Atay and J. Jost, Consensus and synchronization in discrete-time networks of
multi-agents with stochastically switching topologies and time
delays. {\it Networks and Heterogeous Media, 6 (2011), pp. 329-349}.

\bibitem{Liu09}
B. Liu, W. Lu and T. Chen, Consensus in networks of multiagents with
switching topologies modeled as adapted stochastic processes.
{\it SIAM J. Control optim. 49:1 (2011), pp. 227-253}


\bibitem{Hajnal}
J. Hajnal, On products of nonnegative matrices.
{\it Math. Proc. Camb. Phil. Soc. 52 (1976), pp. 521-530.}

\bibitem{Haj2}
J. Hajnal, Weak ergodicity in non-homogeneous Markov
chains. {\it Proc. Camb. Phil. Soc., 54(1958), pp. 233-246}.

\bibitem{Shen}
J. Shen, A geometric approach to ergodic non-homogeneous markov chains.
{\it Wavelet Anal. Multi. Meth., LNPAM, 212 (2000), 341-366. }

\bibitem{Wolfowitz}
J. Wolfowitz, Products of indecomposable, aperiodic, stochastic matrices.
{\it Proceedings of the American Mathematical Society, Vol. 14, No. 5 (Oct., 1963), pp. 733-737.}

\bibitem{Olfati04}
R. Olfati-Saber and R. M. Murray, Consensus problems in networks of agents with switching topology and time-delays.
{\it IEEE Trans. Autom. Control, vol. 49, pp. 1520-1533, 2004.}
\bibitem{Liu10}
X. Liu, W. Lu and T. Chen, Consensus of multi- agent systems with unbounded time-varying delays.
{\it IEEE Transactions on Automatic Control, 55(2010), pp. 2396-2401.}

\bibitem{Liu108}
X. Liu and T. Chen, Consensus problems in networks of agents under nonlinear protocols with directed interaction topology.
 {\it Physical Letters A, 373 (2009), pp. 3122-3127}


\bibitem{schnitzler}
A. Schnitzler and J. Gross, Normal and pathological oscillatory communication in the brain.
{\it Nat. Rev. Neurosci. 6 (2005), pp. 285-296.}

\bibitem{Passino}
K. M. Passino, Biomimicry of bacterial foraging for distributed optimization and control.
{\it IEEE Control Syst. Mag. 22:3(2002), pp. 52-67.}
\bibitem{Montbrio}
E. Montbrio, J. Kurths and B. Blasius, Synchronization of two interacting populations of oscillators.
{\it Phys. Rev. E, 70(2004), 056125.}
\bibitem{Rulkov}
N. F. Rulkov, Images of synchronized chaos: Experiments with circuits.
{\it Chaos, 6:3(1996), pp. 262-279.}
\bibitem{Stone}
L. Stone, R. Olinky, B. Blasius, A. Huppert, and B. Cazelles, Complex Synchronization Phenomena in Ecological Systems.
{\it AIP Conf. Proc.,  622(2002), pp. 476-488.}
\bibitem{Hwang}
K. Hwang, S. Tan, and C. Chen, Cooperative strategy based on adaptive Q-learning for robot soccer systems.
{\it IEEE Trans. Fuzzy Syst., 12(2004), pp.569-.}

\bibitem{Chen11'}
Y. Chen, J. H. Lv, F. L. Han and X. H. Yu, On the cluster consensus of discrete-time multi-agent systems .
{\it Systems and control letters 60(2011), pp. 517-523.}
\bibitem{Yu09}
J. Yu and L. Wang, Group consensus in multi-agent systems with undirected communication graphs.
{\it Proc. 7th Asian Control Conf., 2009, pp. 105-110.}
\bibitem{Yu09'}
J. Yu and L. Wang, Group consensus in multi-agent systems with switching topologies.
{\it IEEE Conference on Decision and Control. 2009, pp. 2652-2657.}
\bibitem{LuChaos}
W. L. Lu, B. Liu. T. P. Chen, Cluster synchronization in networks of coupled nonidentical dynamical
systems, {\it CHAOS, 20 (2010), 013120}.

\bibitem{God}
C. Godsil and G. Royle, ``Algebraic Graph Theory,'' Springer-Verlag, New York, 2001.

\bibitem{Horn}
 R. A. Horn and C. R. Johnson, ``Matrix analysis,'' Cambridge
University Press, 1985.

\bibitem{Wu}
 C. W. Wu, Synchronization and convergence of linear dynamics in random directed networks,
{\it  IEEE Trans. Autom. Control, 51 (2006), pp. 1207-1210}.



\bibitem{Lin}
C.-T. Lin, Structural controlability. {\it IEEE Trans. Auto. Control, 19 (1974), pp. 201-208.}

\bibitem{Rein}
K. J. Reinschke, ``Multivariable Control: A Graph-Theories Approach'', Springer-Verlag, 1988

\bibitem{Dion}
J. -M. Dion, C. Commault, J. van der Woude, Generic properties and control of linear structured systems. {\it Automatica, 39 (2003), pp. 1125-1144.}

\bibitem{Pesin}
 L. Barreira and Y. Pesin, ``Lyapunov exponents and smooth ergodic theory''
 University Lecture Series, AMS, Providence, RI, 2001.

\bibitem{Lu..}
 W. Lu, F. M. Atay and J. Jost, Synchronization of discrete-time networks with time-varying couplings.
 {\it SIAM J. Math. Analys., 39 (2007), 1231-1259.}


\end{thebibliography}
\end{document}